\documentclass{amsart}
\usepackage{amscd}
\newcommand{\R}{\mathbb{R}}
\newcommand{\C}{\mathbb{C}}
\newcommand{\Z}{\mathbb{Z}}
\newcommand{\Q}{\mathbb{Q}}

\newtheorem{theorem}{Theorem}[section]

\theoremstyle{definition}
\newtheorem{definition}[theorem]{Definition}
\newtheorem{example}[theorem]{Example}

\newtheorem{conjecture}[theorem]{Conjecture}

\theoremstyle{remark}
\newtheorem{remark}[theorem]{Remark}

\numberwithin{equation}{section}

\begin{document}

\title{Floer homology of Lagrangian submanifolds}

\author{Kenji Fukaya}
\address{Department of Mathematics, Graduate school of Science, Kyoto University, 
Kitasirahawa, Sakyo-ku, Kyoto 606-8502 Japan}
%\curraddr{Department of Mathematics and Statistics, 
%Case Western Reserve University, Cleveland, Ohio 43403}
\email{fukaya@math.kyoto-u.ac.jp}
\thanks{The author was supported in part by JSPS Grant-in-Aid for Scientific Research
No.18104001 and by Global COE Program G08.}

%\author{Author Two}
%\address{Mathematical Research Section, School of Mathematical Sciences, 
%Australian National University, Canberra ACT 2601, Australia}

%    General info

\date{2011 June 15}



\begin{abstract}
This paper is a survey of Floer theory, which the author 
has been studying jointly with Y.-G. Oh, H. Ohta, K. Ono.
A general idea of the construction is outlined.
We also discuss its relation to (homological) mirror symmetry.
Especially we describe various conjectures 
on (homological) mirror symmetry and various partial results 
towards those conjectures.
\end{abstract}

\maketitle

\section{Introduction}
This article is a survey of Lagrangian Floer theory, which the author 
has been studying jointly with Y.-G. Oh, H. Ohta, K. Ono.
A major part of our study was published as \cite{FOOO08}.
\par 
Floer homology is invented by A. Floer in 1980's. 
There are two areas where Floer homology appears.
One is symplectic geometry and the other is topology of 3-4 dimensional
manifolds (more specifically the gauge theory).
In each of the two areas, there are several different Floer type theories.
In symplectic geometry, there are Floer homology of periodic Hamiltonian system
(\cite{Flo89}) and Floer homology of Lagrangian submanifolds.
Moreover there are two different kinds of Floer theories of contact 
manifolds (\cite{EGH00, Tau07}).
In gauge theory, there are three kinds of Floer homologies: one based 
on Yang-Mills theory (\cite{Flo88II}), 
one based 
on Seiberg-Witten theory (\cite{MaWa01, KrMr}), 
and Heegard Floer homology (\cite{OzSz04}).
Those three are closely related to each other.
\par
All the Floer type theories have common feature that
they define some kinds of   homology theory of ${\infty}/{2}$ degree
 in $\infty$ dimensional space, based on Morse theory.
\par
There are many interesting topics to discuss on the general feature of 
Floer type theories.
I however do not discuss them in this article and 
concentrate on the points which are important 
in Lagrangian Floer theory.

\section{Floer homology of Lagrangian submanifolds}

A {\it symplectic manifold} 
is a pair $(X, \omega)$ of ($2n$)-dimensional manifold $X$ 
and its closed two form $\omega$ such that $\omega^n$ is nowhere $0$ in $X$.
A {\it Lagrangian submanifold} $L$ of $(X, \omega)$ is an $n$-dimensional submanifold  
such that the restriction of 
$\omega$ to $L$ is $0$.
\par
Typical examples of 
$X$ are cotangenet bundle $T^*M$ of a manifold $M$, and a K\"ahler manifold.
Typical examples of $L$ are zero section of $T^*M$, 
and the set of real points of a projective algebraic variety $X$ defined over $\R$.
\par
The Floer homology $HF(L_1, L_2)$ of Lagrangian submanifolds 
is an invariant of a pair $(L_1, L_2)$ of Lagrangian submanifolds 
in a symplectic manifold  $X$.
This is actually the first one \cite{Flo88I} among various  Floer type theories 
that was studied by A. Floer.
However there are several difficulties to establish it in the general 
situation and so it takes much time for such theory to be established.
\par
The ideal properties that Floer homology of Lagrangian submanifolds 
are expected to enjoy can be summarized as follows.

\medskip
\begin{enumerate}
\item[(i)] 
We can associate a module $HF(L_1, L_2)$
to a symplectic manifold $(X, \omega)$ and 
a pair $(L_1, L_2)$ of Lagrangian submanifolds of it.
$HF(L_1, L_2)$ is called the {\it Floer cohomology}.
\item [(ii)]
A pair of Hamitonian diffeomorphims\footnote{The notion of 
Hamitonian diffeomorphim is defined as follows. 
Let $H : [0, 1] \times X \to \R$ be a smooth funciton. We put
$H_t(x) = H(t, x)$. A time dependent Hamiltonian vector field 
 associated to it is a vector $V_{H_t}$ that satisfies  
 $\omega(V, V_{H_t}) = dH_t(V)$, for any 
 vector field $V$. By the non-degeneracy of symplectic form $\omega$
the vector field $V_{H_t}$ is determined uniquely
by this condition. We define 
$\varphi^H_t : X \to X$ by
$\varphi^H_0(x) = x$, 
$
(d/dt)(\varphi^H_0(x))\vert_{t=t_0} = 
V_{H_{t_0}}(\varphi^H_{t_0}(x))
. $
A Hamiltonian diffeomorphism  $\varphi$ is a diffeomorpism 
such that $\varphi = \varphi^H_1$ for some $H$.
It is well known that Hamiltonian diffeomorphism preserves a 
symplectic form.
} $\varphi_i : X \to X$ induces an isomorphism
$HF(L_1, L_2) \cong HF(\varphi_1(L_1), \varphi_2(L_2))$.
\item[(iii)]
If $L_1=L_2=L$ then
$HF(L, L) \cong \bigoplus_{i=0}^n H^i(L)$.
Here the right hand side is the singular homology of $L$.
\item [(iv)]
If $L_1$ intersects transversaly to $L_2$,  then $HF(L_1, L_2)$ 
is generated by at most $\#(L_1\cap L_2)$ elements. Here 
$\#(L_1\cap L_2)$ is the order of the intersection $L_1\cap L_2$.
\end{enumerate}
\medskip
\par
If we assume (i)(ii)(iii)(iv) above,  then
for any pair $(\varphi_1, \varphi_2)$ of diffeomorphisms 
we have 
\begin{equation}\label{eq:lagarnold}
\#(\varphi_1(L) \cap \varphi_2(L)) \ge \sum_{i=0}^n \text{\rm rank}\, H^i(L)
\end{equation}
provided that $\varphi_1(L)$ is transversal to $\varphi_2(L)$.
This is a `Lagrangian version of Arnold's conjecture' 
\footnote{Note however 
that this claim is not correct in general. Arnold certainly did NOT conjecture it
of course.},  
that implies a similar Arnold's conjecture for periodic Hamiltonian system.
\par
Floer extablished $HF(L_1, L_2)$ satisfying 
(i)(ii)(iii)(iv) above, under the assumption,  $\pi_2(X, L) = 0$ and that
there exists Hamiltonian diffeomprphisms $\varphi_i$ 
such that $L_i = \varphi_i(L)$. This assumption is rather restrictive.
\par
Y.-G. Oh \cite{Oh93} relaxed this condition to 
`$L$ is monotone\footnote{
The monotonicity is the condition that 
the two homomorphisms $H_2(M, L;\Z)\to \R$ :
$\beta \mapsto \int_{\beta}\omega$ (where $\omega$ is 
the symplectic form), and the Maslov index $\mu$ (that is a kind of relative version 
of Chern number) are proportional to each other.}
and the minimal Maslov number\footnote{The minimal 
Maslov number is 
the Maslov index of the nonzero holomorphic 
disc with smallest Maslov index.} is not smaller than $3$.'
\par
Actually the Floer homology that satisfies all of (i)(ii)(iii)(iv) above 
can not exist. In fact  for any compact Lagrangian submanifold $L$ 
of $\C^n$ we can find a Hamiltonian diffeomorphism $\varphi$ such that
$L \cap \varphi(L) = \emptyset$. Hence (\ref{eq:lagarnold})
can not hold. ($\C^n$ is non compact. However we can 
take 
$X=\C P^n$ instead. In fact the above mentioned Hamiltonian 
diffeomorphism has a compact support and so is extended to $\C P^n$.
Thus we still have a counter example to (\ref{eq:lagarnold}).)
\par
Our main result on Lagrangian Floer homology 
which modify  (i)(ii)(iii)(iv) is as follows\footnote{In Theorem \ref{thm;mainexistence} 
existence of spin structure is assumed. We may relax it to the 
existence of relative spin structure.}.
Items (1)(2)(3) correspond to (ii)(iii)(iv), respectively.
\begin{theorem}\label{thm;mainexistence}
{\rm (\cite{FOOO00, FOOO08})}
For each spin Lagrangian submanifold $L$ we can associate 
a set $\mathcal M(L)$.
For each pair of spin Lagrangian submanifolds
$L_1, L_2$ and $b_i \in \mathcal M(L_i)$ we can 
associate Floer cohomology $HF((L_1, b_1), (L_2, b_2);\Lambda)$.
They have the following properties.
\smallskip
\begin{enumerate}
\item 
A symplectic diffeomorphism $\varphi : X \to X$ induces 
a bijection $\varphi_* : \mathcal M(L) \to \mathcal M(\varphi(L))$.
A pair of Hamiltonian diffeomorphisms $\varphi_i : X \to X$ $(i=1, 2)$ 
induces the following isomorphism.
\begin{equation}
\aligned
(\varphi_1, \varphi_2)_* : 
&HF((L_1, b_1), (L_2, b_2);\Lambda) \\
&\to HF((\varphi_1(L_1), \varphi_{1*}(b_1)), 
(\varphi_2(L_2), \varphi_{2*}(b_2));\Lambda).
\endaligned
\end{equation}
\item
If $L_1=L_2=L$ we have a spectral sequence $E$ such that
$E_2 \cong H(L;\Lambda)$ and 
\begin{equation}
E_{\infty}^{p} \cong F^p(HF((L, b), (L, b);\Lambda))/
F^{p-1}(HF((L, b), (L, b);\Lambda))
\end{equation}
for an appropriate filtration $F^*(HF((L, b), (L, b);\Lambda))$.
\item
If $L_1$ intersects $L_2$ transversaly,  then 
the rank of $HF(L_1, L_2)$ over $\Lambda$ is not 
greater than $\#(L_1\cap L_2)$.
\end{enumerate}
\end{theorem}

The coefficient ring $\Lambda$ is the {\it universal Novikov field} (\cite{Nov81})
and is the totality of the following infinite series:
\begin{equation}\label{eq;elenov}
\sum_{i=0}^{\infty} a_i T^{\lambda_i}.
\end{equation}
Here $a_i$ are rational numbers and
$\lambda_i$ are real numbers. We assume $\lambda_i$ 
is strictly increasing with respect to $i$ and $\lim_{i\to\infty} \lambda_i = \infty$.
\par
If we assume moreover $\lambda_i$ are rational numbers then 
the totality of such a series is puiseux series ring\footnote{When we 
replace the condition $a_i \in \Q$ by $a_i \in \C$ then 
its becomes the algebraic closure of the formal power series 
ring over $\C$.} and so is a field.
We can define a $T$- adic non Archimedean norm on $\Lambda$ 
and $\Lambda$ is complete with respect to this norm\footnote{In \cite{FOOO08}, we introduced one more 
formal variable $e$ and use a graded ring  $\Lambda_{nov}$.
In this article we do not use $e$. 
In that case Floer homology is only $\Z_2$ graded.}.
\par
The set $\mathcal M(L)$ may be empty.
In that case Theorem \ref{thm;mainexistence}
does not contain any interesting informaiton.
In other words the statement of Theorem \ref{thm;mainexistence}
itself would be rather obvious. (We may simply define  $\mathcal M(L)$ 
is the empty set always. Then all the claims clearly hold.)
So to obtain some nontrivial consequence from Theorem \ref{thm;mainexistence} 
we need to find a condition for 
$\mathcal M(L)$ to be nonempty, or the spectral sequence (2) to degenerate.
We next describe such results. Those results are easier to state after we slightly generalize 
Theorem \ref{thm;mainexistence}.
\par
The set 
$\mathcal M(L)$ is a subset 
$\mathcal M_{\text{\rm weak,def}}(L)$\footnote{The set 
$\mathcal M(L)$ is the set of the gauge equivalence classes of the 
solutions of (\ref{eq;mc}).
We relax the equation (\ref{eq;mc}) 
to $\sum \mathfrak m_k(b^k) \equiv 0 \mod [L]$ 
and $\mathcal M_{\text{\rm weak}}(L)$ is the set of gauge equivalence classes 
of its solutions.
The $A_{\infty}$ structure of 
$H(L;\Lambda_0)$ is deformed by an element of $H^{even}(X;\Lambda_+)$
 (the map (\ref{openGW})),  the equation (\ref{eq;mc})
(or  $\sum \mathfrak m_k(b^k) \equiv 0 \mod [L]$ 
)is deformed accordingly.
The union of the gauge equivalence classes of their solutions is 
$\mathcal M_{\text{\rm weak,def}}(L)$.},  and is described by the map
$$
\pi : \mathcal M_{\text{\rm weak,def}}(L)
\to H(X;\Lambda_+), 
\quad 
\frak{PO} :
\mathcal M_{\text{\rm weak,def}}(L) 
\to \Lambda_+
$$
as,  $\mathcal M(L) = \pi^{-1}(0) \cap \frak{PO}^{-1}(0)$.
The function $\frak{PO}$ is called the {\it potential function}.
Here $\Lambda_0$ is the set of formal sums
(\ref{eq;elenov}) such that 
$\lambda_i \ge0$. It becomes a ring.
$\Lambda_+$ is its maximal ideal  and consists of elements such that $\lambda_i >0$.
The Floer cohomology $HF((L_1, b_1), (L_2, b_2);\Lambda)$ is 
defined 
under the condition for $b_i \in \mathcal M_{\text{\rm weak,def}}(L_i)$ that
$\pi(b_1) = \pi(b_2)$,  $\frak{PO}(b_1) = \frak{PO}(b_2)$.
It satisfies  (i)(ii)(iii)\footnote{(2) holds in case $L_1=L_2$, $b_1=b_2$.}.
\par
\begin{theorem}\label{thm;nonempty}
If the inclusion induced map 
$H^*(X;\Q) \to H^*(L;\Q)$ is surjective for all even $*$, 
then $\mathcal M_{\text{\rm weak,def}}(L)$
is nonempty.
\end{theorem}
We consider the case $L_1=L_2=L$, $b_1=b_2$ in the next theorem.
\begin{theorem}\label{thm;nonvanish}
\begin{enumerate}
\item
The image of the differential of the spectral sequence in $(1)$ Theorem 
$\ref{thm;mainexistence}$
is contained in the Poincar\'e dual to the 
kernel of the map $H_*(L) \to H_*(X)$ induced by the inclusion.
Especially under the assumption of Theorem $\ref{thm;nonempty}$, 
the spectral sequence of Theorem $\ref{thm;mainexistence}$ degenerates
\item 
If the Maslov index of all the discs $(D^2, \partial D^2) \to (M, L)$ is nonnegative 
then the fundamental class of $L$ becomes an nonzero element of $E_{\infty}$
In particular the Floer cohomology $HF((L, b), (L, b);\Lambda)$ is nonzero 
\end{enumerate}
\end{theorem}
The assumptions of the above theorems are rather restrictive.
The condition for the set $\mathcal M(L)$
to be nonempty or the Floer cohomology to be nonzero, 
is closely related to the symplectic topology of $L$.
So it is hard to describe it in term of the topology of $L$ 
only.
Theorem \ref{thm;mainexistence} can be regarded as a 
background result to study the symplectic topology of $L$ using 
Floer cohomology.

\section{Example : toric manifold}
The calculation of Floer cohomology is a difficult problem.
Especially the Lagrangian Floer cohomology is hard to calculate and 
we had only sporadic calculations of it, for a long time.
Recently a systematic calculation becomes possible in the case of toric manifolds.
Let us explain an example before explaining the proof of Theorem \ref{thm;mainexistence}.
\par
A {\it toric manifold} $X$ is a real $2n$ dimensional manifold on which 
real $n$-dimensional torus $T^n$ acts,  so that it admits a map
$\pi : X \to \R^n$ that is called the {\it moment map}.
The fibers of $\pi$ are the $T^n$ orbits,  and the image of $\pi$ is a convex polygon in
$\R^n$. 
We put $\pi(X) = P$. For each $u \in \text{\rm Int}\, P$ the fiber $\pi^{-1}(p)
= L(u)$ is diffeomorphic to an $n$-dimesional torus.
Its Floer  cohomology $HF((L(u), b), (L(u), b);\Lambda)$
is calculated.
In this section we explain a part of this result.
In the case of $L=L(u)$,  the set  $\mathcal M_{\text{\rm weak}}(L(u))$
(that is an element $b$ of $\mathcal M_{\text{\rm weak,def}}(L(u))$
 so that $\pi(b) = 0$)
contains $H^1(L(u);\Lambda_0)$ (See \cite{FOOO08I}).
We restrict $b$ to this subset. We then have a map
$\frak{PO} : H^1(L(u);\Lambda_0) \to \Lambda_0$\footnote{In the last section 
we used $\Lambda_+$ in place of $\Lambda_0$.
This is related to the fact that in case of $\R$ coefficient and $M$ is a real toric 
manifold,  
we can slightly generalized the definition of $\mathcal M_{\text{\rm weak}}(L)$
in \cite{FOOO08I} using the idea of \cite{cho07}.
This generalization is actually possible for any $X$
(\cite{Fuk09}). Since we discuss only in case of $\Lambda_+$
in \cite{FOOO08}, 
we only stated the result over $\Lambda_+$coefficient.}.
\par
In the case of toric manifold,  the potential function $\frak{PO}$ is defined as follows.
We fix a basis $(\text{\bf e}_i)$ of $H^1(L(u);\Z)$, 
and put 
$
\partial \beta = \sum_{i=1}^n (\partial_i \beta)  \text{\bf e}_i
$
for $\beta \in H_2(X, L(u);\Z)$. (Here $\partial_i \beta$ is an integer.)
We describe an element $b$ of $H^1(L(u);\Lambda_0)$
by using the dual basis $(\text{\bf e}^*_i)$ to $(\text{\bf e}_i)$ as 
$b= \sum_{i=1}^n x_i \text{\bf e}^*_i$.
Then, 
\begin{equation}\label{potdef}
\frak{PO}(b) = \frak{PO}(x_1, \ldots, x_n) 
= \sum_{\beta} T^{\beta\cap \omega}\exp\left(\sum_{i=1}^n x_i\partial_i \beta\right) n(\beta).
\end{equation}
Here,  $n(\beta)\in \Q$ is defined roughly as follows.
We fix a point $p\in L(u)$.
Then  $n(\beta)$ is the number of pseudo-holomorphic map 
$(D^2, \partial D^2) \to (X, L)$ that contains  
$p$,  and is of homology class $\beta$.
\begin{theorem}\label{thm;toric}
If  the gradient vector of $\frak{PO}$ is zero at $b \in H^1(L(u);\Lambda_0)$, 
then the Floer cohomology is isomorphic to the singular cohomology.
Namely:
$$
HF((L(u);b), (L(u);b));\Lambda) \cong H(T^n;\Lambda).
$$
Otherwise,  
$$
HF((L(u);b), (L(u);b));\Lambda)=0.
$$
\end{theorem}
The potential function $\frak{PO}$
is closely related to the Landau-Ginzburg super potential 
(See \cite{HKKPTVVZ}). In the case of toric manifold,  
it is mostly calculated in \cite{FOOO08I} based on the result of \cite{cho-oh}.
For example in case $X=\C P^n$ we have $P = \pi(X) 
= \{ (u_1, \cdots, u_n) \mid u_i \ge 0,  \sum u_i \le 1\}$
and the potential function of the fiber $L(u)$ at $u = (u_1, \cdots, u_n)$
is given by:
\begin{equation}\label{eq;potpn}
\frak{PO}(x) = \sum_{i=1}^n T^{u_i} e^{x_i} + T^{1-\sum u_i}e^{-\sum x_i}
\end{equation}
where $x= (x_1, \cdots, x_n) \in \Lambda_0^n \cong H^1(L(u);\Lambda_0)$. 
\begin{example}
We consider the case $S^2 = \C P^1= \{(x, y, z) \mid x^2+y^2+z^2 = 1\}$.
(We normalize its symplectic form so that its area is $1$.)
We consider a $T^1$-action that consists of  rotations around $z$-axis. 
Its orbit is parametrized by $z$-coordinate $z_0$.
We may choose the moment map so that the coordinate $u$ of
$P= [0, 1]$ is the area of 
$\{(x, y, z)\in S^2 \mid z \le z_0\}$.
The complement 
$\C P^1\setminus L(u)$
is divided into two discs,  one $D^2_1(u)$
contains the south pole,  the other 
$D^2_2(u)$ contains the north pole.
We denote the homology class of 
those discs by $\beta_1$,  $\beta_2$ respectively.
These two $\beta$'s are only the discs which contribute 
to the right hand 
side of
(\ref{potdef}). Note
$\partial_1 \beta_1 = 1$,  $\partial_1\beta_2 = -1$, 
and $\beta_1\cap\omega = u$,  $\beta_2 \cap\omega = 1-u$.
Therefore we have
$$
\frak{PO}(x) = T^u e^x + T^{1-u} e^{-x}.
$$
In the case of $\C P^n$,  there are $n+1$ homology classes which 
contribute to the potential function.
\end{example}
The zero of the gradient vector field of the 
potential function $\frak{PO}$ in (\ref{eq;potpn}) exists only in case of 
$u= (1/(n+1), \cdots, 1/(n+1))$. In that case we  have
$x = (\chi, \cdots, \chi)$,  
$\chi = 2\pi\sqrt{-1} (1/(n+1) + \text{\rm integer})$.
\par
The fiber $L(u)$ of this $u = (u_1, \cdots, u_n)$ is called the Clifford torus.
It is well known that for all the other fibers $L(u)$ than Clifford torus,  there exists 
a Hamiltonian diffeomorphism $\varphi$ such that
$\varphi(L(u)) \cap L(u) = \emptyset$.
Therefore 
the Floer cohomology 
$HF((L(u), b), (L(u), b);\Lambda)$ must vanish by
Theorem \ref{thm;mainexistence} (3).
Theorem \ref{thm;toric} is consistent to this fact.
\par\medskip
There are several other calculation of the Lagrangian Floer 
cohomology than toric fibers.
For example in \cite{FOOO09II},  we studied how 
the moduli space of pseudo-holomorphic discs changes by the 
Lagrangian surgery and perform some calculations using it.

\section{$A_{\infty}$ structure}
Floer cohomology has a ring structure. Moreover it is an $A_{\infty}$ algebra.
The set $\mathcal M(L)$ in Theorem \ref{thm;mainexistence} can be 
regarded as a formal scheme defined by those structures.
We need to use $A_{\infty}$ algebra for the proof of Theorem 
\ref{thm;mainexistence} also.
In this section we discuss $A_{\infty}$ structure in Floer theory.
We restrict ourselves to the case of $A_{\infty}$ algebra
associated to a single Lagrangian submanifold.
\par
A $\Z_2$ graded complete 
$\Lambda_0$ module $C$ is said to be a {\it filtered $A_{\infty}$
algebra},  if there exists a series of degree one continuous 
$\Lambda_0$ homomorphisms
$$
\frak m_k : C[1] \otimes_{\Lambda_0} \cdots \otimes_{\Lambda_0}  C[1] \to C[1], 
$$
for $k=0, 1, \cdots$. 
(Here the left hand side is the tensor product of $k$ copies of 
$C[1]$.) such that the relation 
\begin{equation}\label{eq;ainftyrel}
\sum_{k=1}^{n-1}\sum_{i=1}^{k-i+1} (-1)^{*} \frak m_{n-k+1}( 
x_1 \otimes \cdots \otimes \frak m_k(x_i, \cdots, x_{i+k-1}) 
\otimes \cdots \otimes x_n) = 0, 
\end{equation}
holds. 
($*= \deg x_1 + \cdots + \deg x_{i-1} + i - 1$.)
(More precisely we assume some complementary 
conditions such as $\frak m_0(1) \equiv 0 \mod \Lambda_+$.)
Here $C[1]$ is a parity change of $C$,  that is $C[1]^1 = C^0$ and $C[1]^0 =C^1$.
(Note $C$ is $\Z_2$ graded.)
(\ref{eq;ainftyrel}) is called the {\it $A_{\infty}$ relation}.
\par
The operator $\frak m_k$ defines a coderivation by:
$$
\widehat{\frak m}_k : 
x_1 \otimes \cdots \otimes x_n \mapsto
\pm x_1 \otimes \cdots \otimes \frak m_k(x_i, \cdots, x_{i+k-1}) 
\otimes \cdots \otimes x_n.
$$
We put:
$
\hat d = \sum_{k=0}^{\infty} \widehat{\frak m}_k.
$
Then (\ref{eq;ainftyrel}) is equivalent to the formula
$
\hat d\circ \hat d = 0.
$
\begin{theorem}\label{thm:Ainfty}
For each (relatively) spin Lagrangian submanifold $L$, 
its singular cohomology $H(L;\Lambda_0)$ over $\Lambda_0$
has a structre of filtered $A_{\infty}$
algebra.
This filtered $A_{\infty}$ algebra is invariant of the 
symplectic diffeomorphism type of $(M, \omega, L)$
up to isomorphism.
\end{theorem}
We next describe a relation between 
Theorem \ref{thm:Ainfty} and Theorem \ref{thm;mainexistence}.
Let $b \in H^{odd}(L;\Lambda_+)$.
We consider the equation\footnote{The right hand side of 
this equation is infinite sum. However since we assumed $b\equiv 0 \mod \Lambda_+$
it converges in $T$-adic topology.}.
\begin{equation}\label{eq;mc}
\frak m_0(1) + \sum_{k=1}^{\infty} \frak m_k(b^{\otimes k}) = 0.
\end{equation}
Let us consider the case $\frak m_k=0$ for $k\ne 1, 2$. We put
$\frak m_1 = \pm d$,  $\frak m_2(x, y) =\pm x\wedge y$. 
Then (\ref{eq;mc}) becomes the Maurer-Cartan equation $db + b \wedge b = 0$.
 Maurer-Cartan equation gives a condition for the connection $\nabla = d+b$ to be flat.
(Namely it is equivalent to ${\nabla}\circ {\nabla}=0$.)
Therefore  Maurer-Cartan equation 
describes the moduli space of flat bundles\footnote{In (\ref{eq;mc}) the element $b$ is a cohomology classes.
In Maurer-Cartan equation $b$ is a differential form rather than de Rham cohomology 
class.
See three lines right before Remark \ref{remoriryaku} about this point.}.
\par
The set of gauge equivalence classes\footnote{We omit the definition of gauge equivalence. 
See \cite{FOOO08} section 4.3.} of (\ref{eq;mc}) is the set $\mathcal M(L)$.
For $[b_1], [b_2] \in \mathcal M(L)$,  $x \in H(L;\Lambda_0)$ we define
$$
\delta_{b_1, b_2}(x) = \sum_{k_1, k_2} \frak m_{k_1+k_2+1}
(b_1^{\otimes k_1}\otimes x \otimes b_2^{\otimes k_2}). 
$$
Then (\ref{eq;ainftyrel}), (\ref{eq;mc}) implies
$\delta_{b_1, b_2}\circ \delta_{b_1, b_2} = 0$.
The cohomology group of this boundary operator $\delta_{b_1, b_2}$
is Floer cohomology $HF((L, b_1), (L, b_2);\Lambda_0)$\footnote{We use $\Lambda$ instead of $\Lambda_0$
in Theorem \ref{thm;mainexistence}.
This is because (2) holds only over this coefficient ring.}.
\par
To prove Theorem \ref{thm:Ainfty} we define operators
$\frak m_k$ on the cohomology group. The definition is 
based on the moduli space of pseudo-holomorphic discs and 
is performed  as follows.
\par
We first take an almost complex structure $J : TX \to TX$ on $X$ 
that is an $\R$ linear map such that $J\circ J = -1$.
We assume that $J$ is compatible with $\omega$,  namely 
$g(V, W) = \omega(V, JW)$ is a Riemannian metric.
We say a map $u : \Sigma \to X$
is {\it pseudo-holomorphic} or $J$-holomorphic if the formula:
$$
J(du(V)) = du(j(V))
$$
is satisfied.
Here $j :  T\Sigma \to T\Sigma$ is a complex structure of
$\Sigma$.
Let $\beta\in H_2(X, L;\Z)$. For
$k\ge 0$ we define the set
$\widetilde{\mathcal M}_{k+1}^{\text{reg}}(L;\beta)$
to be the totality of all pairs\footnote{$D^2 = 
\{z \in \C \mid \vert z\vert \le 1\}$.} $(u, \vec z)$ of maps $u : D^2 \to X$ and
$\vec z=(z_0, \cdots, z_k)$ that have the following 6 properties.
\begin{enumerate}
\item $u$ is pseudo-holomorphic.
\item $u(\partial D^2) \subset L$. 
\item The relative homology class of $u$ is $\beta$.
\item $z_i \in \partial D^2$. 
\item If $i\ne j$ then $z_i \ne z_j$.
\item $z_0, z_1, \cdots, z_k$ respects the counterclockwise 
cyclic order of $\partial D^2$.
\end{enumerate}
\par
The group of all biholomorphic isomorphisms of $D^2$
is isomorphic to $PSL(2;\R)$.
It acts on $\widetilde{\mathcal M}_{k+1}^{\text{reg}}(L;\beta)$ by
$g\cdot(u, \vec z) = (g\cdot u, g\vec z)$, 
$(g\cdot u)(z) = u(g^{-1}(z))$, 
$g\vec z = (gz_0, \cdots, gz_k)$.
\par
We denote by ${\mathcal M}_{k+1}^{\text{reg}}(L;\beta)$
the quotient space of this action.
We define a quotient topology on it.
It has a compactification which is 
Hausdorff.
To obtain a compactification, we use a generalization of the 
notion of the stable map (\cite{Kon95II}) to the case of bordered Riemann surface. 
We denote the compactification by
${\mathcal M}_{k+1}(L;\beta)$.
The map
$
ev = (ev_0, \cdots, ev_k) : 
\widetilde{\mathcal M}_{k+1}^{\text{reg}}(L;\beta) 
\to L^{k+1}
$
is defined by $ev_i(u, \vec z) = u(z_i)$.
This map is invariant of $PSL(2;\R)$ action and extends to the compactification.
In other words we have a map
$
ev : {\mathcal M}_{k+1}(L;\beta) \to L^{k+1}.
$
We call it the {\it evaluation map}. 
Let us assume the following for simplicity.
\par\medskip
\noindent (*) \quad ${\mathcal M}_{k+1}(L;\beta)$ is a compact oriented manifold.
\par\medskip
We remark that this assumption is rarely satisfied.
We discuss this points in Sections 5-6.
In case (*) is satisfied, 
the dimension of the moduli space ${\mathcal M}_{k+1}(L;\beta)$
is given by:
\begin{equation}\label{eq;dimformula}
\dim {\mathcal M}_{k+1}(L;\beta) = n + (k+1) + \mu(\beta) - 3
= n+ k - 2 + \mu(\beta).
\end{equation}
Here $\mu : H_2(X, L;\Z) \to \Z$ is a homomorphism called Maslov index.
\par
In case 
(*) is satisfied,  we define
$
\frak m_{k, \beta} : (H^*(L;\Z))^{k\otimes} \to H^*(L;\Z)
$
by
\begin{equation}\label{eq;defm}
\aligned
&\frak m_{k, \beta}(P_1, \cdots, P_k) \\
&= PD_L(ev_{0*}(PD_{{\mathcal M}_{k+1}(L;\beta)}((ev_1, \cdots, ev_k)^*(P_1\times \cdots \times P_k)))).
\endaligned
\end{equation}
Here $PD$ is the Poincar\'e duality map. The first one is $PD_L : H_d(L) \to H^{n-d}(L)$, 
and the second one is
$$
PD_{{\mathcal M}_{k+1}(L;\beta)} : H^{d}({\mathcal M}_{k+1}(L;\beta))
\to H_{n+ k - 2 + \mu(\beta)-d}({\mathcal M}_{k+1}(L;\beta)).
$$
\par
Then we define $\frak m_k$ by
\begin{equation}\label{defmk}
\frak m_k = \sum_{\beta} T^{\beta\cap \omega}\frak m_{k, \beta}.
\end{equation}
Here $\beta\cap \omega = \int_{u}\omega$.
($u$ is a relative cycle representing $\beta \in H_2(X, L;\Z)$.)
Since $\omega$ is $0$ on $L$,  the integral $\int_{u}\omega$
depends only on the relative homology class 
$\beta$ and is independent of $u$.
We can use 
Gromov compactness (that is the compactness of the moduli space of pseudo-holomorphic maps $u$ so that
$\int_{u}\omega$ is smaller than a fixed number) to show 
the convergence of $\frak m_k$ with respect to the non-Archimedean topology of $\Lambda_{0}$.

\section{Kuranishi structure}
In the last section, we assumed (*) to define $\frak m_{k, \beta}$.
The assumption (*) is rarely satsified. For example it is never satisfied in case $\beta=0$, $k>2$.
As a consequence we {\it cannot} use (\ref{eq;defm}) itself to define
$\frak m_{k, \beta}$. 
We use the method of virtual fundamental chain\footnote{The notion of virtual fundamental class is a natural generalization 
of the notion of the fundamental class $[M] \in H_n(M;\Z)$ of an 
oriented manifold $M$.
In our situation, we do not take the homology class 
but use cycle or chain. So we need the notion of virtual fundamental {\it chain}.}
and chain level intersection theory to overcome this trouble.
We discuss the former in this section and the later in the next section.
We need to study the following two points.
\par\smallskip
\noindent(A) In general, an element of ${\mathcal M}_{k+1}(L;\beta)$ has a singularity.
Moreover it may have nontrivial automorphism. 
\par
\noindent(B) The moduli space ${\mathcal M}_{k+1}(L;\beta)$
has a (codimension 1) boundary. So 
even in the situation (A) does not occur, it does not determine a cycle.
\par\smallskip
Let us first discuss the point (A).
In order to isolate this point from (B), 
we study the following situation,  which is slightly 
different from one in the last section.
We consider a (real) 6 dimensional symplectic manifold 
$X$ with $c_1(X)=0$.
(We use a compatible almost complex structure $J$ 
to define the first Chern class $c_1(X)$.)
For each $\alpha \in H_2(X;\Z)$,
we consider the set of the pseudo-holomorphic maps 
$u : S^2 \to X$ of homology class $\alpha$.
We identify the maps which are transformed by 
the group of automorphisms of $S^2$($= PSL(2;\C)$) 
and compactify the set of 
equivalence classes of this identification. Then we obtain the moduli space 
$\mathcal M(\alpha)$. In case 
$\mathcal M(\alpha)$ is a manifold, we calculate 
its dimension using 
$c_1(X) = 0$ and obtain
$\dim \mathcal M(\alpha) = 0$.
In other words,  $\mathcal M(\alpha)$ is an oriented $0$ dimensional compact manifold.
So its order (counted with sign) makes sense. This is the simplest case of Gromov-Witten 
invariant $GW(\alpha)$.
\par
The problem (A) in this case is the problem to 
justify `the order counted with sign'  when
$\mathcal M(\alpha)$ may be not a manifold.
\par
The theory of virtual fundamental class resolves this problem.
In the case of $GW(\alpha)$ there 
are two methods to work it out. 
One is based on algebraic geometry (\cite{LiTi96,  BF})
and the other is based on differential geometry (\cite{FuOn99II, Rua99, LiTi98, Sie96}).
In the case when Lagrangian submanifold is included, 
we can not use the method of algebraic geometry.
So we explain the method based on differential geometry.
\par
We use the notion of Kuranishi structure to 
define virtual fundamental class.
A Kuranishi structure on a topological space $\mathcal M$
is,  roughly speaking,  a system to give the way to represent 
the space  $\mathcal M$ locally as the zero set of an equation 
(on a finite dimensional space) so that they are glued in a consistent way.
Namely for each  $u\in\mathcal M$, 
we represent a neighborhood of $u$ in $\mathcal M$ as the set of 
the solutions of $m$ equations
$s_{u, i}(x_1, \cdots, x_n)=0$ ($i=1, \cdots, m$)  of $n$-variables.
The main idea is to include not only the solution set but also the equation itself as a part of the structure.
When we change the base point $u$ the equation may change. 
However only the following 2 kinds of changes are allowed. 
\begin{enumerate}
\item We increase the number of equations and variables by the same 
number,  say $\ell$. The equation itself is modified by the following 
trivial way.
$$\aligned
s_{u, i}(x_1, \cdots, x_n, x_{n+1}, \cdots, x_{n+\ell}) &= s_{u, i}(x_1, \cdots, x_n) 
\qquad i=1, \cdots, m \\
s_{u, i+m}(x_1, \cdots, x_n, x_{n+1}, \cdots, x_{n+\ell}) &= x_{i+n}
\qquad i=1, \cdots, \ell
\endaligned$$
\par
\item The coordinate transformation of 
the variables $x_i$ and the linear transformation of the equations.
We may allow the linear transformations 
of the equation to depend on $x_i$.
Namely we allow the following transformation $(s_u;x_1, \cdots, x_n) \mapsto (s_v;y_1\cdots, y_n)$.
$$
s_{v, j}(y_1, \cdots, y_n)
= \sum_{i=1}^m g_{ji}(x_1, \cdots, x_n) s_{u, i}(x_1, \cdots, x_n),  \quad
y_i 
= y_i(x_1, \cdots, x_n).
$$
\item We include the process to divide the whole structure by a finite group.
\end{enumerate}
\par\medskip
We include the map $s_u$ in the data defining Kuranishi structure. 
So the multiplicity is determined from the Kuranishi structure.
For example,   in case $n=m=2$, $s(x, y) = (x^2-y^2, 2xy)$, 
the solution set $s^{-1}(0)$ consists of one point $0$.
This is the same as the case $s(x, y) = (x, y)$.
However they are different as Kuranishi structures.
The virtual fundamental class (which we define later) is a $0$
dimensional homology class,  that is a number,  in this case. In case $s(x, y) = (x^2-y^2, 2xy)$, it is $2$.
In case $s(x, y) = (x, y)$,  it is $1$.
\par\medskip
Let us explain the notion of Kuranishi structure in more detail 
and also explain the way how it appears in the study of pseudo-holomorphic curve.
\par
When we study moduli spaces using differential geometry, 
we study the set of the solutions of a differential equation, 
which is elliptic modulo the action of the `gauge group'.
Here we study the nonlinear Cauchy-Riemann equation:
\begin{equation}\label{eq;nlcr}
J \circ du = du \circ j.
\end{equation}
It linearization is a differential operator:
$$
D_u\overline{\partial} : \Gamma(S^2;u^*TX) \to \Gamma(S^2;u^*TX 
\otimes \Lambda^{01}(S^2)), 
$$
where $\Lambda^{01}(S^2)$ is the complex 
line bundle of $(0, 1)$ forms on $S^2$, 
and $\Gamma$ denotes the set of smooth sections.
\par
If $D_u\overline{\partial}$ is surjective 
then we can use implicit function theorem 
to show that the solution set of (\ref{eq;nlcr}) 
is a smooth manifold in a neighborhood of $u$.
It is finite dimensional by ellipticity.
We denote a neighborhood of $u$ in the solution set by $V(u)$.
\par
Our moduli space is the quotient space of the set of solutions of (\ref{eq;nlcr}) by 
the $PSL(2;\C)$-action. If the map
$u$ is nontrivial,  the set of $g \in PSL(2;\C)$ satisfying $g\cdot u = u$ is 
finite. We denote it by $\Gamma(u)$.
Then the neighborhood of $[u]$ in our moduli space $\mathcal M(\alpha)$
is 
$V(u)/\Gamma(u)$. The space
$V(u)$ is a manifold if $D_u\overline{\partial}$ is surjective.
Therefore $V(u)/\Gamma(u)$ is a quotient of a manifold by a finite group.
Satake introduced the notion of a space which is locally a quotient 
of a manifold by a finite group.
It is nowadays called an {\it orbifold}.
Compact and oriented orbifold without boundary 
carries a fundamental class in the same way as manifold.
Thus Problem (A) does not occur if  $\mathcal M(\alpha)$
is an orbifold.
\par
The problem (A) occurs when $D_u\overline{\partial}$ is not surjective.
The local theory of such a moduli space is studied by Kuranishi 
in the case of the moduli space of complex manifolds.
We apply it in the case of the moduli space 
of pseudo-holomorphic curves.
\par
We take a finite dimensional linear subset $E(u) \subset \Gamma(S^2;u^*TX 
\otimes \Lambda^{01}(S^2))
$ such that the sum of the image of $D_u\overline{\partial}$ and $E(u)$
generates
$\Gamma(S^2;u^*TX 
\otimes \Lambda^{01}(S^2))$ as  a vector space.
This is possible since $D_u\overline{\partial}$ is a Fredholm operator.
\par
We may choose a family of isomorphisms 
$\Gamma(S^2;u^*TX 
\otimes \Lambda^{01}(S^2)) \cong
\Gamma(S^2;v^*TX 
\otimes \Lambda^{01}(S^2))$ depending smoothly on $v$.
Then we may regard $E(u)$ as a subset of $\Gamma(S^2;v^*TX 
\otimes \Lambda^{01}(S^2))$.
Now we replace the equation (\ref{eq;nlcr})  by
\begin{equation}\label{eq;nrcrmod}
J\circ dv - dv \circ j \in E(u).
\end{equation}
We can then apply implicit function theorem to it.
We regard the element of $\Gamma(S^2;v^*TX 
\otimes \Lambda^{01}(S^2))$ as a map which 
associate to
$x \in S^2$ an anti-complex-linear map $T_{x}S^2 \to T_{v(x)}X$.
We can make sense (\ref{eq;nrcrmod}) by this identification.
\par
Therefore the set of solutions of (\ref{eq;nrcrmod}) is 
a finite dimensional manifold,  which we denote by $V(u)$.
We define a map $s_u : V(u) \to  E(u)$ by
$
v \mapsto J\circ dv - dv \circ j
$.
Finally we put
$$
\Gamma(u) = \{ g \in PSL(2;\C) \mid u\circ g = u\}.
$$
(This is the automorphism group.)
The group $\Gamma(u)$ is a finite group.
\par
We take $E(u)$ so that it is $\Gamma(u)$ invariant. 
Then $\Gamma(u)$ acts on $V(u)$ and $E(u)$,  so that 
$s_u$ is $\Gamma(u)$-invariant.
We thus obtain the following:
\par
A neighborhood $V(u)$ in a complex vector space of finite dimension.
A linear and effective action of  a finite group $\Gamma(u)$ on it.
A finite dimensional representation of $\Gamma(u)$ on 
$E(u)$.
Moreover we have a $\Gamma(u)$-equivariant map $s_u : V(u) \to  E(u)$
and a homeomorphism 
$
\psi_u : s_u^{-1}(0)/\Gamma(u) \to \mathcal M(\alpha)
$
onto a neighborhood of $u$.
\par
We call 
$(V(u), E(u), \Gamma(u), s_u, \psi_u)$ a {\it Kuranishi chart}.
\par
We can define the notion of coordinate change between 
Kuranishi charts in the same way as the definition of manifold 
as follows.
\par
A {\it coordinate change} from a Kuranishi chart
$(V(u_1), E(u_1), \Gamma(u_1), s_{u_1}, \psi_{u_1})$ to another Kuranishi chart
$(V(u_2), E(u_2), \Gamma(u_2), s_{u_2}, \psi_{u_2})$ consists of 
a $\Gamma(u_1)$-invariant open subset $V(u_1, u_2) \subset V(u_1)$, 
a group homomorphism $\phi_{u_2u_1} : \Gamma(u_1) \to \Gamma(u_2)$, 
$\phi_{u_2u_1}$-equivariant smooth embedding $\varphi_{u_2u_1} : V(u_1, u_2) \to V(u_2)$, 
a $\phi_{u_2u_1}$-equivariant smooth embedding of vector bundles
$\hat\varphi_{u_2u_1} : V(u_1, u_2)\times E(u_1) \to V(u_2) \times E(u_2)$
over $\varphi_{u_2u_1}$ such that
the following compatibility condition is satisfied:
$$
s_{u_2} \circ \varphi_{u_2u_1} = \hat\varphi_{u_2u_1} \circ s_{u_1}, 
\qquad
\psi_{u_2} \circ \varphi_{u_2u_1} = \psi_{u_1} \quad \text{ holds on $s_{u_1}^{-1}(0)\cap V(u_1, u_2)$ 
}.
$$
We can define the compatibility conditions 
between coordinate changes in the same way as the 
definition of manifolds.
\par
A paracompact Hausdorff space 
is said to have a {\it Kuranishi structure}
when it is covered by 
Kuranishi charts so that 
coordinate changes are defined among them 
and the compatibility between them is satisfied.
We assume one more condition (\ref{virdim}) explained later.
We call $E(u)$ the {\it obstruction bundle},  and $s_u$ the 
{\it Kuranishi map}.
\par
The biggest difference between manifold structure and Kuranishi structure 
lies in the fact that we do not assume the embedding $\varphi_{u_2u_1} : V(u_1, u_2) \to V(u_2)$
to be a local homeomorphism.
Especialy the dimension of $V(u_1)$ may not be equal to the dimension of $V(u_2)$.
We assume the following condition instead. 
\begin{equation}\label{virdim}
\text{$\dim V(u) - \mathrm{rank}\,  E(u)$ is independent of $u$.}
\end{equation}
We call this difference the {\it dimension} (or 
{\it virtual dimension}) 
of our Kuranishi structure.
\par
We need one more condition to define the 
virtual fundamental class of the Kuranishi structure.
We consider the normal bundle $N_{V(u_1, u_2)}V(u_2)$
of our smooth embedding $\varphi_{u_2u_1} : V(u_1, u_2) \to V(u_2)$.
The compatibility condition implies that $s_{u_2}$
induces the following linear map between vector bundles.
$$
ds_{u_2} : N_{V(u_1, u_2)}V(u_2)=\frac{\varphi_{u_2u_1}^*TV(u_2)}{TV(u_1)} \to \frac{\varphi_{u_2u_1}^*E(u_2)}{E(u_1)}.
$$
We say that our Kuranishi structure has a {\it tangent bundle} if 
$ds_{u_2}$ is an isomorphism.
From now on we consider only the Kuranishi structure 
with tangent bundle.
We can define the notion of orientability or orientation
of Kuranishi structure with tangent bundle.
\par
As we mentioned before, a Kuranishi structure on 
$\mathcal M$ may be regarded as a system 
which assigns a way to represent $\mathcal M$ as the solution 
set of the equation
$
s_{u} = 0
$
on a neighborhood of each point $u$,  so that 
it is consistent when we move  $u$.
The notion of a coordinate of a space with singularity is 
classical in algebraic geometry.
The definition of algebraic variety or scheme is bases on such a notion\footnote{Since we 
include the finite group action our notion of Kuranishi structure 
corresponds to Deligne-Mumford stack in algebraic geometry.}.
In complex analytic category, 
there is a similar notion such as analytic space or 
Douady space.
Here we work on $C^{\infty}$ category.
As a consequence the following point is different from those.
\par
Scheme is define as a locally ringed space. There is a difficulty to do so here.
\par
For each Kuranishi chart $(V(u), E(u), \Gamma(u), s_u, \psi_u)$ we put
$s_u = (s_u^1, \cdots, s_u^m)$. Then we may 
consider the quotient 
of the ring of function germs $C^{\infty}_0(V(u))$
by the ideal $(s_u^1, \cdots, s_u^m)$ generated by 
$s_u^1, \cdots, s_u^m$ and define 
a local ring $C^{\infty}_0(V(u))/(s_u^1, \cdots, s_u^m)$.
By moving 
$u$ we obtain a locally ringed space.
However this locally ringed space is a bit hard to study.
The biggest problem is that $C^{\infty}_0(V(u))$
is harder to handle compared to the 
ring of germs of holomorphic functions or 
polynomial rings. 
For example the Krull dimension of $C^{\infty}_0(V(u))$ is infinite.
In particular the Krull dimension  of $C^{\infty}_0(V(u))/(s_u^1, \cdots, s_u^m)$
is not equal to 
$\dim V(u) - \mathrm{rank}\,  E(u)$.
\par
By this reason we avoid using the ring $C^{\infty}(V(u))/(s_u^1, \cdots, s_u^m)$
or its localization. 
Instead we cover our space by the charts which have positive size. 
Moreover we assumed compatibility among the Kuranishi maps.
In the situation when we can use the language of locally ringed space, 
the compatibility of the Kuranishi maps can be replaced by the 
condition that the structure sheaf is defined globally.

\begin{remark}
Here is another difference between our case and the case of
algebraic geometry or analytic space.
Suppse $\dim V = 1$,  $\operatorname{rank} E = 1$.
We consider the following three cases.
`$s(x) = x$',  `$s'(x) = x ^3$',  `$s''(x) = e^{-1/x} (x>0),  -e^{1/x} (x<0)$'.
The multiplicity of $s$ and $s'$ at origin is $1$ and $3$ respectively 
in the usual sense of algebraic geometry.
So we might define the virtual fundamental cycle of 
$(\R, \R, \{1\}, s)$ to be $1$,  and of $(\R, \R, \{1\}, s')$ to be $3$.
However we define the virtual fundamental cycle of $(\R, \R, \{1\}, s')$ 
to be $1$.
This is justified by considering $s'_{\epsilon}(x) = x^3 + \epsilon x$ for 
sufficiently small $\epsilon$.
(Note $x$ is a real variable.)
\par
In the case of $(\R, \R, \{1\}, s'')$,  the multiplicity 
might be infinity if we consider the analogy of algebraic geometry.
However the virtual fundamental cycle $(\R, \R, \{1\}, s'')$ is $1$.
\end{remark}
Let us define the virtual fundamental cycle of the space $\mathcal M$
with oriented Kuranishi structure.
A map $f: \mathcal M \to Y$ to a topological space $Y$ is said 
to be 
{\it strongly continuous} if $f$ is extended to its 
Kuranishi neighborhood $f_{u} : V(u) \to Y$ for each
$u \in \mathcal M$ and 
they are compatible with $\varphi_{u_2u_1}$.
We omit the detail of the compatibility condition.
The virtual fundamental class is an element of $H_{\dim \mathcal M}(Y)$.
\par
We first consider the case when $\Gamma(u)$ is always trivial.
Then we can deform the Kuranishi map $s_u$ on each of 
the Kuranishi chart $(V(u), E(u), \Gamma(u), s_u, \psi_u)$
to obtain $\frak s_u$,  such that 
$\frak s_u^{-1}(0)$ becomes a $\dim \mathcal M$ dimensional manifold.
Using the compatibility of the Kuranishi charts we can take $\frak s_u$ 
so that they are compatible with the coordinate change in a suitable sense.
Then the zero sets,  $\frak s_u^{-1}(0)$ are glued to define a $\dim \mathcal M$
dimensional manifold. 
The strong continuity implies that $f$ defines a map from this manifold to $Y$.
Therefore 
$f_*[\cup_u \frak s_u^{-1}(0)] \in H_{\dim \mathcal M}(Y;\Z)$ is defined.
\par
In case 
$\Gamma(u)$ is nontrivial,  it is impossible to find a $\Gamma(u)$-invariant 
$\frak s_u$ that is transversal to $0$, in general.
In this case, we use multisection instead.
We first take $\frak s_u$ which is transversal to $0$.
Then we consider the totality of the 
$\Gamma(u)$ orbits of it. Namely we consider $\{\gamma \cdot \frak s_u \mid \gamma \in \Gamma(u) \}$.
This is $\Gamma(u)$ invariant as a set. We regard it as a multivalued section.
Then its zero set, that is the set of all points 
where at least one of the $\#\Gamma(u)$-branches of this multivalued section is zero, 
carries a fundamental cycle as follows.
We first triangulate this zero set. 
Then for each simplex $\Delta_a$ of top dimension,  we define a wight $m_a 
\in \Q$ as follows. We divide the number of the 
branches which becomes zero on it by the order $\#\Gamma(u)$. The weight is 
this ratio. $\sum m_a(\Delta_a, f)$ is a singular chain of $Y$.
We can show that it is a cycle of $Y$.
This is the {\it virtual fundamental class} $\in H(Y;\Q)$.
\par
In case $\dim\mathcal M=0$ we do not need to specify $Y$ to define 
the virtual fundamental class as a degree $0$ homology class.
It is a rational number. We regard it as the `number of the points of 
$\mathcal M$'.
\par
The moduli space $\mathcal M(\alpha)$ of pseudo-holomorphic curves
has an oriented Kuranishi structure.
In case 
$[u] \in \mathcal M(\alpha)$ is represented by a map $u : S^2 \to X$ we 
explained the way to find its Kuranishi neighborhood already.
We need a compactification and the case $u : S^2 \to X$ 
corresponds to the case when $u$ is in the interior of 
$\mathcal M(\alpha)$.
To define a Kuranishi structure on $\mathcal M(\alpha)$
we need to define a Kuranishi chart on a neighborhood of the 
point
$u$ at infinity.
Such $u$ is a map from a singular Riemann surface $\Sigma$ with 
normal crossing singularity.
Starting with $u : \Sigma \to X$, we can define  a 
family of pseudo-holomorphic maps from the normalization 
of $\Sigma$. This  is an important topic in the theory 
of pseudo-holomorphic curve and
is called the gluing.
(A similar procedure had been studied in gauge theory.
It was initiated by Taubes to prove an existence theorem
of self-dual connection on 4 manifolds.)
We use the theory of glueing to construct a Kuranishi chart on a neighborhood of a 
point
$u$ at infinity.
\begin{remark}\label{varioussimiar}
The notion of Kuranishi structure was introduced at the year 1996
by Fukaya-Ono. The construction of virtual 
fundamental class by the differential geometric method was done 
independently by Li-Tian,  Ruan,  Siebert in the same year.
Fukaya-Ono used the notion of multisection explained above.
Ruan used de Rham cohmology and Li-Tian used the notion of normal 
cone which had been used in algebraic geometry in a related but different 
purpose.
Fukaya-Ono and Ruan used finite dimensional approximation.
Li-Tian did not take finite dimensional approximation and studied infinite 
dimensional space directly.
\par
After the turn of the century, 
a notion called polyfold is proposed by Hofer-Wysocki-Zehnder.
The theory of polyfold follows ours 
in the two points below.
\begin{enumerate}
\item It defines an appropriate class of spaces 
including various moduli spaces. 
It associates a virtual fundamental cycle 
to the spaces in that class in a way independent 
of the way how such a structure is obtained.
\item
It uses multivalued abstract perturbation\footnote{
Here abstract perturbation is the method in which, 
instead of specifying the explicit way to perturb the equation, 
we consider the moduli space locally as a zero set 
of an abstract map and perturb it in an abstract way.
This is not so much a new idea and actually was used in 
\cite{Don} more than 30 years ago.}.
\end{enumerate}
It is also easy to show the following:
For any polyfold there is a space with Kuranishi structure which has 
the same virtual fundamental class.
\par
Therefore the story of Kuranishi structure 
is applicable to {\it any} problem to which the story of polyfold is applicable.
\par
It seems that the most novel part of the theory of polyfold is 
its analytic part.
The theory of Kuranishi structure starts 
at the borderline where analysis is over and topology starts.
Namely in the story of Kuranishi structure the construction
of Kuranishi structure is left to the study of each of the 
problems and the abstract theory starts at the point 
where finite dimensional approximation (the Kuranishi structure) 
is constructed.
On the other hand,  in the story of polyfold, one of the 
main part of the construction that is the glueing is included in
the general theory and it formulates the situation where glueing becomes possible.
Since the story of polyfold is not yet worked out in detail,  we do not discuss it here.
\end{remark}
\begin{remark}\label{kurakate}
As we mentioned before, the notion of Kuranishi structure 
is a kind of $C^{\infty}$ analogue of the notion of scheme and stack.
Therefore it seems important to define a category of the space with 
Kuranishi structure\footnote{
From the point of view of analogy with scheme,  
fiber product is important. The fiber product 
of spaces of Kuranishi structure over manifold is defined in \cite{FOOO08} \S A1
and is applied. We remark however that this may not be the fiber product in the sense of 
category theory}.
Especially it is important to find a good notion of morphisms between them.
The author did not find a good way to do so yet.
A difficulty to do so is as follows.
A coordinate change of Kuranishi chart should be an example 
of such morphisms and should be an isomorphism.
Namely we need to regard the Kuranishi structure $(U, E, s)$ and 
$(U \times \R^m,  E \times \R^m,  (s, id))$ to be isomorphic.
We can define a map $(U, E, s)\to (U \times \R^m,  E \times \R^m,  (s, id))$ 
in a natural way. However it is hard to find a map in the opposite direction 
in the usual sense. So we need to localize the category. 
Several conditions are required for the localization to be well-defined. 
Those conditions are not trivial to check.
\end{remark}
\section{Chain level intersection theory}
The story described in the last section is the case 
of pseudo-holomorphic map from closed Riemann surface without boundary, 
that was establishd in 1996.
The theory of Floer homology 
studies the case of the 
pseudo-holomorphic map from compact Riemann surface with boundary
(bordered Riemann surface). 
The novel point appearing here is problem 
(B) in the last section. We discuss it in this section
(based on \cite{FOOO08} section 7.2).
\par
The main problem is that the moduli space $\mathcal M_{k+1}(L;\beta)$ has
Kuranishi structure with corners.
The notion of Kuranishi structure with corners is 
defined in the same way as before by including the case when 
$U(u)$ is an open subset of $[0, \infty)^k \times \R^{n-k}$.
\par
Note a manifold with boundary or corner does not carry a 
fundamental class. In other words the correspondence by 
a manifold with boundary or corner  
does not induce a map between homology groups.
By this reason, we need chain level intersection theory.
\begin{remark}
If a strongly continuous map $f :  \mathcal M \to Y$ and 
a subset $Z \subset Y$ satisfy
$f(\partial \mathcal M) \subset Z$, then we can define a 
`relative virtual fundamental cycle' $f_*([\mathcal M]) \in H(Y, Z;\Q)$ in the same way.
However in Lagrangian Floer theory, we need to consider the case when
$ev : \mathcal M_{k+1}(L;\beta) \to L^{k+1}$ is the strongly 
continuous map $f$. In this case,  there does not seem to be a 
reasonable proper subset of $L^{k+1}$ containing $ev(\partial(\mathcal M_{k+1}(L;\beta)))$.
The author does not know any example 
where `relative virtual fundamental cycle' in the above sense was applied successfully.
\end{remark}
Because of the problem we mentioned above, 
we need to perform the construction of the 
operator $\frak m_k$ in the chain level,  in order to construct our 
$A_{\infty}$ algebra. 
The construction is roughly as follows.
For each  singular chain $\sigma_i : \Delta^{d_i} \to L$, 
we take the fiber product 
\begin{equation}\label{fiberdelta}
\mathcal M_{k+1}(L;\beta) \, \, \, {}_{ev_1, \cdots, ev_k}\times_{\sigma_1, \cdots, \sigma_k} 
\, \,  \prod_{i=1}^k\Delta^{d_i} 
\end{equation}
over $L^{k}$ and triangulate it. Then we regard each of the simplices of the top dimension of 
(\ref{fiberdelta}) as a singular chain by the map $ev_0$.
The sum of them is $\frak m_{k, \beta}(\sigma_1, \cdots, \sigma_k)$ by definition.
It is rather a heavy job to work out its detail\footnote{The longest section \S 7.2 of \cite{FOOO08} is devoted to work it out.}.
We explain this construction a bit more below.
\par
We denote the fiber product 
(\ref{fiberdelta}) as $\mathcal M_{k+1}(L;\beta;\sigma_1, \cdots, \sigma_k)$.
This space has a Kuranishi structure with corner.
Its boundary is decomposed into the sum of the following 
two types of spaces:
\begin{enumerate}
\item [(a)]
The fiber product between 
$
\Delta^{d_1} \times \cdots \times \Delta^{d_{i-1}} \times   
\mathcal M_{j-i+1}(L;\beta_1;\sigma_i, \cdots, \sigma_j)
\times \Delta^{d_{j+1}} \times \cdots \times \Delta^{d_k}
$
and
$\mathcal M_{k+1+i-j}(L;\beta_2)$
\item [(b)]
$
\mathcal M_{k+1}(L;\beta;\sigma_1, \cdots, \partial_c\sigma_i, \cdots, \sigma_k)
$
Here $\partial_c\sigma_i$ is $c$-th face of the sigular chain $\sigma_i$.
\end{enumerate}
Then we define multi-sections on $\mathcal M_{k+1}(L;\beta;\sigma_1, \cdots, \sigma_k)$
and triangulations of its zero set
inductively according to the order of 
$k$, $\omega\cap\beta$, $\sum(n-d_i)$, 
so that they are compatible with the above description of the 
boundaries. It determines the virtual fundamental chains of (\ref{fiberdelta}) for each of $\beta, \sigma_1, \cdots, \sigma_k$. We thus defined
$
\frak m_{k, \beta}(\sigma_1, \cdots, \sigma_k)
$
as a singular chain. We define $\frak m_{1, 0}$ to be the boundary operator of the 
singular chain complex. Then $\frak m_k$ is defined by (\ref{defmk}).
\par
We can prove that $\frak m_k$ defines an $A_{\infty}$ structure as follows.
The boundary of (\ref{fiberdelta}) is 
$$
\frak m_{1, 0}(\frak m_{k, \beta}(\sigma_1, \cdots, \sigma_k)).
$$
This is decomposed into (a) and (b) above.
The case (a) gives 
$$
\frak m_{k+i-j, \beta_1}(\sigma_1, \cdots, \sigma_{i-1}, \frak m_{j-i+1, \beta_2}(\sigma_i, \cdots, \sigma_{j}), \sigma_{j+1}, \cdots, \sigma_{k})
$$
and (b) gives 
$$
\frak m_{k, \beta}(\sigma_1, \cdots, 
\frak m_{1, 0}(\sigma_i), \cdots, \sigma_k).
$$
Thus we have
$$\aligned
&\frak m_{1, 0}(\frak m_{k, \beta}(\sigma_1, \cdots, \sigma_k))\\
&= \sum_i \pm\frak m_{k, \beta}(\sigma_1, \cdots, \frak m_{1, 0}(\sigma_i), \cdots, \sigma_k) \\
&\quad+\sum_{\beta_1+\beta_2=\beta}\sum_{1\le i\le j \le k} \pm\frak m_{k+i-j, \beta_1}(\sigma_1, \cdots, \sigma_{i-1}, \frak m_{j-i+1, \beta_2}(\sigma_i, \cdots, \sigma_{j}), \sigma_{j+1}, \cdots, \sigma_{k})
\endaligned$$
This is equivalent to (\ref{eq;ainftyrel}).
\par
This construction is nontrivial also when we restrict it to the case of $\beta=0$.
In that case it defines a structure of $A_{\infty}$ algebra on the singular chain complex.
\begin{remark}
We here explained the construction using singular homology.
At the time of writing this article,  two other constructions are known.
\par
One is to use de Rham cohomology.
In this case,  the formula (\ref{fiberdelta}) corresponds to
$$
(\rho_1, \cdots, \rho_k)
\mapsto (ev_0)_!\left(ev_1^*\rho_1 \wedge \cdots \wedge ev_k^*\rho_k\right).
$$ 
Here $\rho_i$ is a differential form on $L$ and $ev_i^*\rho_i$ 
is its pull back. $(ev_0)_!$ is the integration along the fiber 
by the map $ev_0: \mathcal M(L;\beta) \to L$.
Even in the case when 
$\mathcal M(L;\beta)$ is a smooth manifold, the map $ev_0$ may not be a submersion.
So the integration along the fiber defines a  distributional form that may not be 
a smooth form.
So we can not define the operator $\frak m_k$ on the de Rham complex of 
$L$ in this way. We use smoothing of the distributional form $\mathcal M(L;\beta)$
so that it is compatible at the boundary,  inductively.
Smoothing differential forms are performed by using a continuous family 
of multisections. See \cite{FOOO08} Section 7.5 and
\cite{Fuk07I,  Fuk09}.
\par
An advantage of this method is that it is easier to 
keep symmetry. On the other hand,  the author does not know how to 
work over $\Q$ coefficient when we use de Rham cohomology.
\par
The other method is to use the notion of 
Kuranishi homology proposed by \cite{Joy}.
Namely we regard a pair $(\mathcal M, f)$ of the space with
Kuranishi structure $\mathcal M$ and a strongly continuous 
(weakly submersive) map $f:\mathcal M \to Y$ as a chain on $Y$
Joyce called it a Kuranishi chain and construct 
the structure on the chain complex consisting of Kuranishi chains.
We can define the fiber product among Kuranishi chains always.
So using this method we can construct a 
structure on the chain complex of  Kuranishi chains without 
perturbing the equation at all\footnote{Perturbation 
becomes necessary to study the relation between 
Kuranishi homology with other homology theories 
such as singular homology.}.
A difficulty of this construction is as follows:
\par
When we define a homology theory by regarding 
`a pair of a space and a map from it' as a chain,  we need to 
eliminate the automorphism by some method. 
Otherwise it does not give a correct homology group.
\par
Let us consider the pair of a manifold with corner $N$ and a map $f : N \to Y$.
We identify $(N, f) \sim (N', f\circ h)$ for each diffeomorphism 
$h : N' \to N$.
We take the set of equivalence classes of this identifications and 
let $SM(Y)$ be the free abelian group whose basis is identified 
with this set. We define a boundary operator 
$\partial(N, f) = (\partial N, f)$.
We thus obtain a chain complex. However its homology is not
isomorphic to the ordinary homology of $Y$.
(The reason is,  roughly speaking,  the diffeomorphism provides too 
big freedom of identifications.
\par
Joyce resolved this problem by including `gauge fixing data' as a 
part of the data of Kuranishi chain and eliminate the automorphism.
\end{remark}
The construction of this section is not canonical 
and many choices are involved during the construction.
However we can show that the resulting $A_{\infty}$  algebra 
on the cohomology group $H(L;\Lambda_0)$ is 
independent of those choices up to an isomorphisms of 
filtered $A_{\infty}$ algebra.
To prove it we construct an appropriate 
$A_{\infty}$ homomorphism between the $A_{\infty}$ 
structures on singular chain complex and show that 
it is  a homotopy equivalence.
We do so by inductively constructing chain maps using 
an appropriate moduli spaces.
\begin{remark}
To prove the well-defined-ness of the structure up to homotopy equivalence,  
we first need to build a homotopy theory of filtered $A_{\infty}$ algebra.
There are various references (\cite{Lef03, Smi00, Kel01}) describing homotopy theory 
of $A_{\infty}$ algebra (without filtration).
However,  for example, it is hard to find a reference where
the equivalence of various  definitions of homotopy between 
$A_{\infty}$ homomorphisms is proved in detail. 
(The author knows at least two different definitions.)
See \cite{FOOO08} Chapters 4 and 5 for the homotopy theory of filtered $A_{\infty}$ algebra.
Homotopy theory of $A_{\infty}$ algebra is regarded as 
a generalization of the homotopy theory of differential graded algebra 
(\cite{Sul78}).
\end{remark}
Once we obtain a structure of  filtered $A_{\infty}$ algebra on the 
singular chain complex then we can use homological algebra 
to squeeze it to the homology group.
(This is a classical result which goes back to Kadei\v svili \cite{Kad82}. 
See also \cite{KoSo01}.)
We thus obtain a  structure of filtered $A_{\infty}$ algebra on 
the cohomology group.
\begin{remark}\label{remoriryaku}
We omit several important parts of the proof.
Especially we omit the argument on the orientation and sign.
The (relative) spin structure is used to orient the moduli spaces 
$\mathcal M(L;\beta)$.
We remark that proving the orientability of $\mathcal M(L;\beta)$
is only the first (and rather easier) step of the whole works 
on sign and orientation.
To work out the sign and orientation part of the construction of the 
$A_{\infty}$ structure,  we need the following:
Choose orientations 
of many spaces so that their fiber products are related at the 
corners and the boundaries.
To check those orientations are consistent to the sign 
appearing in the homological algebra of 
$A_{\infty}$ structures. This is heavy and cumbersome job,  
and is performed in Chapter 8 \cite{FOOO08}, 
which occupies around 80 pages.
\end{remark}

\section{$A_{\infty}$ category and homological mirror symmetry}
There are two kinds of applications of Lagrangian Floer theory.
One is to the symplectic geometry and the other is 
to the mirror symmery.
Application to the symplectic geometry 
is described for example in \cite{Ono08}.
So we focus here to the application to the mirror symmetry.
A important problem where Lagrangian Floer theory is related to mirror symmetry 
is M. Kontsevitch's homological mirror symmetry conjecture \cite{Kon95I}.
We discuss it in this article.

\begin{definition}\label{Afinitycate}
A {\it filtered $A_{\infty}$ category} $\mathcal C$ consists of 
the set of objects $\mathfrak{Ob}(\mathcal C)$, 
the set of morphisms $\mathcal C(c, c')$ for each $c, c'\in \mathfrak{Ob}(\mathcal C)$, 
such that $\mathcal C(c, c')$ is a $\Lambda_0$ module, and 
a $\Lambda_0$ module homomorphisms
$$
\frak m_k : 
\bigotimes_{i=1}^k \mathcal C(c_{i-1}, c_i) \to \mathcal C(c_{0}, c_k)
$$
for each $c_0, \cdots, c_k\in \mathfrak{Ob}(\mathcal C)$, 
$k=1, \cdots, \infty$.
We assume the relation (\ref{eq;ainftyrel}) among them.
\end{definition}
\begin{example}
In case $\frak m_k=0$ for $k\ne 2$, the $A_{\infty}$ category 
becomes a usual additive category.
Note $\frak m_2$ is different from the composition 
by sign.
\par
The case $\frak m_k=0$ for $k\ne 1, 2$ is called the 
{\it differential graded category}.
It is introduced by
\cite{BK}.
\end{example}
To each symplectic manifold $(X, \omega)$, we can associate an 
filtered $A_{\infty}$ category $\mathcal{LAG}(X)$ whose object
is a pair $(L, b)$ of a spin Lagrangian submanifold
$L$ and $b \in \mathcal M(L)$. In case of 
$c_i = (L, b_i)$, namely in the case when Lagrangian submanifolds 
$L$ are the same, we define 
$\mathcal{LAG}(X)((L, b_i), (L, b_{i+1}))\cong H(L;\Lambda_{0})$ and
$$
\frak m_k(x_1, \cdots, x_k)
= \sum_{\ell_0=0}^{\infty}\cdots\sum_{\ell_k=0}^{\infty}
\frak m_{\ell_0+\cdots+\ell_k + k}
(b_0^{\otimes\ell_0}, x_1, b_1^{\otimes\ell_1}, \cdots, b_{k-1}^{\otimes\ell_{k-1}}, x_k
, b_{k}^{\otimes\ell_{k}}).
$$
Here the right hand side is $\frak m$ in Theorem \ref{thm:Ainfty}, and the left hand 
side is $\frak m$ in Definition \ref{Afinitycate}.
The detail of the definition of this $A_{\infty}$ category is in \cite{Fuk02II}.
\par
$A_{\infty}$ category is not an abelian category.
However we can replace the notion of chain complex in abelian category 
by the twisted complex and can define its derived category \cite{Kon95I}.
\par
Here the twister complex is defined as follows. Let
$c_i\in \mathfrak{Ob}(\mathcal C)$, $i=1, \cdots, n$, 
$x_{ij} \in \mathcal C(c_i, c_j)$, 
$1\le i \le j \le n$.
We say $(\{c_i\}, \{x_{ij}\})$ is a {\it twisted complex} if, for each
$1 \le a, b \le k$, the formula 
\begin{equation}\label{twistedcomplex}
\sum_k\sum_{j_0=a, \cdots, j_k=b}\frak m_k(x_{j_0j_1}, \cdots, x_{j_{k-1}j_k})  = 0
\end{equation}
is satisfied.
\par
In case $\frak m_k = 0$ for $k\ne 2$ and $x_{ij} = 0$ for $j\ne i+1$, 
(\ref{twistedcomplex}) becomes the relation $\frak m_2(x_{i+1 i+2}, x_{i i+1}) = 0$.
Namely it defines a chain complex.
\par
We can define a mapping cone of a morphism of twisted complex.
We thus obtain a triangulated category.
See \cite{Fuk02II, Sei06}. We call this triangulated category the
{\it derived category} of $\mathcal C$.
\par
The homological mirror symmetry conjecture 
by Kontsevitch \cite{Kon95I} is as follows.
\begin{conjecture}\label{CYmirror}
For each Calabi-Yau manifold $X$
we can associate another 
Calabi-Yau manifold $\hat X$, its mirror, 
such that 
the derived category of the category of coherent analytic sheaves on $\hat X$ is equivalent to the 
derived category of $\mathcal{LAG}(X)$.
\end{conjecture}
\begin{remark}
A Calabi-Yau manifold is a K\"ahler manifold and so has 
both symplectic structure (K\"ahler form) and the complex structure.
In Conjecture \ref{CYmirror} we consider the symplectic structure of 
$X$, and the complex structure of $\hat X$, only.
The former is called the {\it A-model} and the later is called the 
{\it B-model} \cite{Witten}.
\end{remark}
The statement of the above conjecture is slightly imprecise.
Let us explain this point first. The problem is the following: 
In the category $\mathcal{LAG}(X)$ the set of morphisms is a module over the 
universal Novikov ring $\Lambda_0$.
On the other hand, the set of morphisms in the category of coherent analytic 
sheaves of 
$X$ is a complex vector space.
So they can not be isomorphic.
There are two ways to correct this point.
\begin{enumerate}
\item 
We substitute a sufficiently small positive number for $T$.
Then the set of morphisms of $\mathcal{LAG}(X)$
becomes a $\C$ vector space.
For this purpose we need to show that the 
structure constants of $\frak m_k$ converges 
when we substitute a sufficiently small positive number for $T$.
This is actually very difficult to prove.
\item
We modify the category of coherent analytic sheaves on $\hat X$ so that 
the set of morphisms becomes a module over the Novikov ring \cite{KoSo01, Fuk03}.
\end{enumerate}
To realize the plan (2), we regard $\hat X$ not as a single Calabi-Yau manifold 
but a family of it parametrized by a disc $D^2(\epsilon)$ with 
small diameter. Namely  we consider a family 
$\pi : \frak X \to D^2(\epsilon)$.
\par
We assume all the fibers of 
$\pi$ other than $\pi^{-1}(0)$ is nonsingular.
The projection 
$\pi$ is not arbitrary.
In mirror symmetry the case when $0$ is a maximal degenerate point 
appears. Here
$0$ is a said to be a maximal degenerate point, 
if $\pi^{-1}(0)$ is a union of irreducible components which are 
normal crossing and that there exists a point in $\pi^{-1}(0)$ where 
$\dim_{\C}X+1$ irreducible components meet.
For example, let us define a family of degree $n+1$ hyper surfaces of $\C P^n$ by
$$
\{ ([x_0:\cdots:x_n], t) 
\in \C P^n \times D^2 \mid x_0\cdots x_n = t(x_0^{n+1}+\cdots + x_n^{n+1}) \}.
$$
Then the fiber of 
$\pi : ([x_0:\cdots:x_n], t) \mapsto t$ at $t$ is nonsingular for $t\ne 0$
and is singular at $t=0$.
Moreover the $n$ irreducible components meet at a point in 
the fiber of  $t=0$.
\par
Now we suppose that 
$0$ is a maximal degenerate point.
We formalize the family along the fiber of $0$.
We then obtain a formal scheme over $\C[[t]]$.
The set of the morphisms of the category of its 
coherent sheaves is a module over $\C[[t]]$.
When we include the branched covering whose 
branch locus is in 
$0$, then the coefficient ring becomes the Puiseux series ring, that is 
very close to the Novikov ring.
Later on in various part of the story 
we may either work on Novikov ring or on $\C$.
We need some nontrivial arguments to go from one to the other.
We however omit the argument to do so\footnote{Actually we do not 
know the way to go from one to the other completely.}.
\par
We can state a part of the homological mirror symmetry 
more explicitly as follows.
\begin{conjecture}\label{CYmirror2}{(\cite{Fuk02III})}
For each pair $(L, b)$ of Lagrangian submanifold $L$ of $X$ and $b\in \mathcal M(L)$, we can associate a chain complex $\mathcal E(L, b)$ of the 
coherent analytic sheaves over $\hat X$, such that the following holds.
\begin{enumerate}
\item There exists an isomorphism\footnote{${\rm Ext}$ is the derived functor of the functor which 
associate $Hom(\mathcal E_1, \mathcal E_2)$ to a pair of the coherent analytic sheaves $\mathcal E_1$, $\mathcal E_2$.}
$$
HF((L_1, b_1), (L_2, b_2)) \cong {\rm Ext}(\mathcal E(L_1, b_1), \mathcal E(L_2, b_2)).
$$
\item The following diagram commutes.
$$
\begin{CD}
HF((L_1, b_1), (L_2, b_2))\otimes HF((L_2, b_2), (L_3, b_3)) @>{\frak m_2}>>
HF((L_1, b_1), (L_3, b_3)) \\
@V{\cong}VV  @V{\cong}VV
\\
{\rm Ext}(\mathcal E(L_1, b_1), \mathcal E(L_2, b_2))\otimes {\rm Ext}(\mathcal E(L_2, b_2), \mathcal E(L_3, b_3)) @>>>
{\rm Ext}(\mathcal E(L_1, b_1), \mathcal E(L_3, b_3))
\end{CD}
$$
Here the lower horizontal arrow is the Yoneda product and the 
vertical arrows are the isomorphisms in (1)\footnote{In Theorem \ref{thm:Ainfty} $\frak m_2$ is defined on the singular cohomology 
$H(L(u);\Lambda_0)$.
Since $\frak m_2$ is a derivation with  respect to the boundary operator $\frak m_1$, 
it follows that $\frak m_2$ defines a product structure on the $\frak m_1$ cohomology, 
that is the Floer cohomology.}.
\end{enumerate}
\end{conjecture}
\section{homological mirror symmetry and classical mirror symmetry}
The classical mirror symmetry
is a statement which claims the coincidence of the generating function 
of the number of pseudo-holomorphic curves on $X$ and 
the generating function obtained from the deformation theory 
of complex structures (Yukawa coupring) of $\hat X$. 
In this section, we explain its relation to homological mirror symmetry.
\par
For each pair of $A_{\infty}$ categories $\mathcal C_1$, $\mathcal C_2$, 
there exists an $A_{\infty}$ category $\mathcal{FUNC}(\mathcal C_1, \mathcal C_2)$
such that its objects is an $A_{\infty}$ functor $\mathcal C_1 \to \mathcal C_2$ \cite{Fuk02II}.
\par
The identity functor $\mathcal C \to \mathcal C$ is an object of 
$\mathcal{FUNC}(\mathcal C, \mathcal C)$, which we denote by $1_{\mathcal C}$.
The set of (pre) natural transformations from $1_{\mathcal C}$ to 
$1_{\mathcal C}$ becomes an $A_{\infty}$ algebra \cite{Fuk02II}. 
We denote it by $\frak{Hom}(1_{\mathcal C}, 1_{\mathcal C})$.
If $\mathcal C$ is an $A_{\infty}$ category with one object only, that is nothing but an
$A_{\infty}$ algebra $C$, then $\frak{Hom}(1_{\mathcal C}, 1_{\mathcal C})$ coincides with the 
Hochschild complex 
$
CH(C, C) = (\bigoplus_k {Hom}(C^{\otimes k}, C), \delta)
$.
(In case $C$ is an associative algebra it coincides with the Hochschild complex  in the 
usual sense.)
In the case of general $A_{\infty}$ category $\mathcal C$, we call
$\frak{Hom}(1_{\mathcal C}, 1_{\mathcal C})$ the 
{\it Hochschild complex} also.
\begin{conjecture}\label{QCtoHoch}
Under certain assumption on 
$X$\footnote{For example the author believes that all the
projective Calabi-Yau manifolds 
have this property.},  the
Hochschild cohomology
$H(\frak{Hom}(1_{\mathcal C}, 1_{\mathcal C}))$ 
of $\mathcal C = \mathcal{LAG}(X)$ is isomorphic to the 
quantum cohomology ring $QH(X;\Lambda)$ of
$X$.
\end{conjecture}
Conjecture \ref{QCtoHoch} claims that $\mathcal{LAG}(X)$ determines 
the quantum cohomology ring.
\par
On the other hand, we can define a
similar Hochschild complex from the derived category of 
the category of coherent analytic sheaves on $\hat X$.
The cohomology of this Hochschild complex and its 
product structure determines the deformation theory of $\hat X$ and 
Yukawa coupling on the deformation space of it.
Thus Conjecture \ref{QCtoHoch} implies the equality 
\par\medskip
\centerline{quantum cup product = Yukawa coupling, } 
\par\medskip
\noindent that is 
the classical mirror symmetry\footnote{
We did not give a precise description how the quantum cohomology 
is related to the deformation theory of its mirror, 
but only suggests that there should be some relation.
The author does not know more precise way to state this relation.
The author does not know a reference where such 
relation is established or at least conjectured in a precise way either.}.
The discussion above seems to be written in \cite{Kon95I}.
Seidel \cite{Sei02} states it as Conjecture 4 including the case when $X$ 
is not necessarily compact.
\par
Let us discuss Conjecture \ref{QCtoHoch} more.
The map
\begin{equation}\label{openGW}
QH(X;\Lambda_0) \to H(\frak{Hom}(1_{\mathcal C}, 1_{\mathcal C}))
\end{equation}
can be defined by Open-Closed Gromov-Witten theory (See \cite{FOOO08} section 3.8.)
as follows.
Let us consider the set $\widetilde{\mathcal M}_{k+1}^{\text{reg}}(L;\beta)$ 
we defined in section 4.
We divide it by the action of 
$
U(1) \cong \{ g \in {\rm Aut}(D^2, j_{D^2}) \cong PSL(2;\R) \mid g(0) = 0\}
$
and compactify it. We denote it by $\mathcal M_{1, k+1}(L;\beta)$.
We can define $ev : \mathcal M_{1, k+1}(L;\beta) \to L^{k+1}$.
Moreover we define 
$$
ev^{\text{int}} : \mathcal M_{1, k+1}(L;\beta) \to M
$$
by 
$$
ev^{\text{int}}(u) = u(0).
$$
If $Q$ is a cycle of $M$ then we can use fiber product to define:
$$
 \mathcal M_{1, k+1}(L;\beta;Q) =  \mathcal M_{1, k+1}(L;\beta)\, \, {}_{ev^{\text{int}}}\times_M Q.
$$
This space has a Kuranishi structure. We can use it 
in place of $\mathcal M_{k+1}(L;\beta)$ and proceed in the same way as section 6.
We then obtain:
$$
\frak q(Q;\cdots) : H(L;\Lambda_0)[1]^{\otimes k} \to H(L;\Lambda_0)[1].
$$
Namely we obtain an element of $CH(HF((L, b), (L, b)), HF((L, b), (L, b)))$ for each of $L, b$.
They behave in a functorial way when we move $L, b$. It thus defines an element of 
$H(\frak{Hom}(1_{\mathcal C}, 1_{\mathcal C}))$.
We associate it to $[Q] \in H(M)$ and obtain the map (\ref{openGW}).
We remark that $QH(X;\Lambda_0)$ is isomorphic to the sigular cohomology as a module over 
$\Lambda_0$. (The ring structure is deformed.)
We can perform the above construction and prove that (\ref{openGW})
is a ring homomorphism without assuming any extra condition on $X$.
The proof of later is similar to the proof of the associativity of 
the quntum cup product.
Seidel proved it in the case Lagrangian submanifold is exact\footnote{Namely there exists one form $\theta$ on $M$ such that 
$\omega=d\theta$ for the symplectic form $\omega$, and there exists a function 
$f$ on $L$ such that $\theta = df$ on $L$.}.
In \cite{BC}\footnote{Note Floer cohomology is called the Lagrangian quantum homology in 
\cite{BC}. I have a strong objection to the 
authors of \cite{BC} who try to chang the name of the notion which had been established 
long time ago
and try to eliminate the name of Floer who discovered this important notion.} the 
case when $L$ is monotone is proved.
\par
The hardest part of the proof of the Conjecture \ref{QCtoHoch} 
is the proof that 
(\ref{openGW}) is an isomorphism.
This does not hold unless we assume some condition to $X$.
In fact, there is a symplectic torus 
which does not contain a Lagrangian submanifold without 
nontrivial Floer cohomology.
For such $X$, the homomorphism (\ref{openGW}) is not an isomorphism.
Proving that 
(\ref{openGW}) is an isomorphism is 
proving the existence of enough many 
Lagrangian submanifolds so that they distinguish all the 
cohomology classes of $X$.
It may be regarded as a `mirror to the Hodge conjecture' and 
is a very difficult problem to solve in general.
There are various cases where the map 
(\ref{openGW}) is proved to be an isomorphism.
It is proved that (\ref{openGW}) is an isomorphism
in the case of toric manifold in \cite{FOOO08X}.
\par
The homological algebra of the map 
(\ref{openGW}) is closely related to a conjecture by Deligne
which is solved in \cite{KoSo00}.

\section{Strominger-Yau-Zaslow conjecture}
We next describe a relation between homological mirror symmetry 
conjecture and a conjecture by 
Strominger-Yau-Zaslow \cite{SYZ96}.
Let us consider a way to construct $\hat X$ so that the homological mirror 
symmetry\footnote{From the point of view of \cite{SYZ96}, 
it may be better to say D-brane duality.} 
 holds.
For each point $p$ of a complex manifold $\hat X$, we can 
define a skyscraper sheaf $\mathcal F_p$.
Namely we put $\mathcal F_p(U)=\C$ if $p\in U$
and $\mathcal F_p(U)=0$ if $p\notin U$.
Let us assume that there exists $(L, b)$ such that 
$\mathcal E(L, b) = \mathcal F_p$ by the correspondence in 
Conjecture \ref{CYmirror2}. We put $(L_p, b_p)=(L, b)$
Then Conjecture \ref{CYmirror2} (1) implies
$$
HF((L_p, b_p), (L_p, b_p)) \cong {\rm Ext}(\mathcal F_p, \mathcal F_p).
$$
We can calculate the ${\rm Ext}$ group of the skyscraper sheaf easily and 
can show that the right hand side is isomorphic to $H(T^n;\C)$ that is the cohomology 
of the $n$-dimensional torus.
Floer cohomology $HF((L_p, b_p), (L_p, b_p))$ is related to the cohomology of 
$L_p$ by the spectral sequence in 
Theorem \ref{thm;mainexistence} (2).
In case the spectral sequence degenerates we have
$H(L_p) \cong H(T^n)$.
Thus we conjecture that the Lagrangian submanifold 
$L_p$ is a torus.
\par
We next consider $b_p$.
In general it is an odd degree cohomology class of $L_p$.
When we study Calabi-Yau manifold, the Floer cohomology is not only $\Z_2$ graded but  is $\Z$ graded
(See \cite{Sei00}). We need to assume $b_p \in H^1(L_p)$ for $\Z$ grading.
In this case $b_p$ corresponds to an element of 
$Hom(H_1(L_p;\Z)), \C\setminus\{0\})$ that is a flat connection.
It is believed that this flat connection is unital\footnote{When we work over $\C$ coefficient (and not over Novikov ring coefficient) 
the parameter which deforms the connection to non-unitary is 
$Hom(H_1(L_p;\Z), \R_{>0})$ and becomes a part of the 
parameter to deform the Lagrangian submanifold
$L$.}.
In sum we have:
\par\medskip
The mirror manifold $\hat X$ is a moduli space 
of the pair of Lagrangian torus $L_p$ and a flat $U(1)$ connection on it.
\par\medskip
Lagrangian torus in a symplectic manifold appears in 
the study of integrable system.
Namely if $(X, \omega)$ is a symplectic manifold and $\pi :  X \to B$ is a
map whose fiber is a compact Lagrangian submanifold, 
then we can show that the fiber of 
$\pi$ is a torus\footnote{
This follows from Liouville-Arnold theorem which asserts that 
if there are $n$-independent first integral and if the 
level set of them is compact, then the level set is a torus.
The coordinate of the fiber is called the angle coordinate 
and the coordinate of $B$ is called the action coordinate.}.
\par
For each $q \in B$ we put $L_q = \pi^{-1}(q)$.
The moduli space of flat $U(1)$ connections on  $L_q$ becomes  
the dual torus of $L_q$.
In other words, the mirror $\hat X$ of $X$ is a union of the dual tori 
of the fibers.
\par
We may sum up the above discussion to the 
following conjecture by  Strominger-Yas-Zaslow.
We need to include the case when there is a singular fiber.
We consider a map $\pi :  X \to B$ from a 
$2n$-dimensional symplectic manifold $(X, \omega)$ to an $n$-dimensional manifold $B$, 
so that the fiber of the generic point is an $n$-dimensional Lagrangian torus.
Let $B_0$ be the subset of $B$ consisting of the points 
whose fiber is an $n$-dimensional Lagrangian submanifold.

\begin{conjecture}\label{SYZ}
For each maximal degenerate family of Calabi-Yau manifolds, we have 
a projection $\pi :  X \to B$ such that 
$B\setminus B_0$ is of codimension $\ge 2$ in $B$.
We obtain a family $\hat X_0 \to B_0$ by taking fiber-wise dual.
$\hat X$ is a compactification of $\hat X_0$.
\end{conjecture}
\begin{remark}
Strominger-Yas-Zaslow conjectures the fibers $L_q$ to be  special 
Lagrangian submanifolds\footnote{Especially it is conjectured to be a minimal submanifold.}.
Nowadays it is known that we need some modification on this points.
See \cite{Joy07, DG09}. 
\end{remark}
In the case when Conjecture \ref{SYZ} holds, 
we expect to obtain the assignment 
$(L, b) \mapsto \mathcal E(L, b)$ in Conjecture \ref{CYmirror2} 
as follows.
For simplicity we assume that $L$ is transversal to the fibers of
$\pi :  X \to B$. A point $p$ of 
$\hat X$ is identified with $(L_{q(p)}, b(p))$ by mirror symmetry. Here ${q(p)} \in B$ and
$L_{q(p)}$ is its fiber. The element 
$b(p)$ determines a flat bundle on it and we may regard it as 
an element of $\mathcal M(L_{q(p)})$.
\begin{conjecture}\label{HMS}(\cite{Fuk02I, Fuk05I})
The coherent sheaf $\mathcal E(L, b)$
is obtained by the holomorphic bundle on $\hat X$ whose fiber at
$p$ is identified with the Floer cohomology  $HF((L_{q(p)}, b(p)), (L, b))$.
\end{conjecture}
Conjecture \ref{HMS} provides a way to construct $\mathcal E(L, b)$, as follows.
We consider the family of 
Floer cohomologies $HF((L_{q(p)}, b(p)), (L, b))$ where $p \in \hat X$ moves.
If we can define a holomorphic structure on it so that it becomes a holomorphic 
vector bundle on 
$\hat X$, then we can {\it define} 
$\mathcal E(L, b)$ to be this holomorphic 
vector bundle. 
\par
This idea was discovered during the authors' discussion with M. Kontsevitch 
at the year 1998 during his stay in IHES.
It was realized 
in the case of Abelian variety or complex torus in 
\cite{Fuk02III}. Note in those cases, there are no 
singular fiber. So the situation is simpler.
We remark that 
in the case of elliptic curve, where we can calculate the 
both sides of the (homological) mirror symmetry conjecture directly, homological 
mirror symmetry was studied earlier by \cite{Kon95I, PoZa98}.
We also refer \cite{KoSo01, AbSm} for homological mirror symmetry of 
complex torus.
\par
Let us discuss the case of general Calabi-Yau manifold 
where there is a sigular fiber.
We describe the way how we obtain a 
complex structure on $\hat X_0$ in Conjecture \ref{SYZ}.
We expect to perform the construction in two steps.
In the first step, we use flat affine structure on $B$ 
and local tensor calculus to define a 
complex structure. (This complex structure is 
called {\it semi-flat}.) 
In the second step, we add the corrections that are induced by the 
pseudo-holomorphic discs\footnote{There is no correction in the 
case of complex torus.}. This correction is 
called the {\it instanton correction}.
\par
The semi-flat complex structure is defined as follows.
We consider $\pi :  X \to B$. For $q \in B_0$ we have 
a canonical isomorphism $H^1(L_q;\R) \cong T_qB$.
(Here $L_q = \pi^{-1}(q)$.) We use the lattice $H^1(L_q;\Z)$ of
$H^1(L_q;\R)$ to define a flat affine structure on $B$.
\par
On the other hand, 
the moduli space of flat unitary connections 
of the fibers is locally the cohomology group $H^1(L_q;\sqrt{-1}\R)$ with 
coefficient in the Lie algebra $\sqrt{-1}\R$ of
$U(1)$.
Conjecture \ref{SYZ} asserts that these two determine $\hat X$
locally. Therefore the tangent space 
$T_q\hat X$ is identified with $H^1(L_q;\R) \oplus H^1(L_q;\sqrt{-1}\R)$, 
that has a complex structure induced by $J_0$
\footnote{There are various reference on the construction of 
this complex structure $J_0$.
The author quote as many reference as he knows on this point 
in section 2 of \cite{Fuk05I}.
So we omit those references here.}.
\par
We next describe the instanton correction.
The existence of instanton correction is related to the singular fiber as follows.
The semi-flat complex structure 
$J_0$ can be constructed only on $\hat X_0$.
Its compactification $\hat X$ contains a (dual to) the singular fibers.
The complex structure $J_0$ however does not extend to $\hat X$.
Therefore we need a correction to 
extend it to a complex structure on $\hat X$.
Existence of such an instanton correction or quantum correction 
had been know in the physics literature. (See for example \cite{MR97i:81117}.)
This phenomenon is related to the wall crossing of the 
Floer cohomology as was 
observed in \cite{Fuk02III, FOOO00}.
One method to study it is to reduce it to Morse homotopy\footnote{In other words, we use the coincidence of the moduli space of pseudo-holomorphic discs 
in the cotangent bundle and the set of the solution of an ordinary equation 
to a map from a graph.},  that is equivalent to the 
enumeration of pseudo-holomorphic curves via tropical geometry.
Relation to Morse homotopy is discussed in \cite{FuOh97, KoSo01, Fuk05I}.
The case of cotangent bundle is discussed in \cite{FuOh97}.
In \cite{KoSo01} it is claimed that one can generalize \cite{FuOh97} to the 
case when the fiber is a torus. 
The study of singular fiber then is started in \cite{Fuk05I}.
See \cite{Mik06} for tropical geometry.
\par
In \cite{Fuk05I} the author discussed the way how the 
instanton  correction is related to the Gromov-Witten invariant.
There Dolbeault cohomology is used to describe 
deformation of complex structure.
In \cite{Fuk05I}, the deformation of $\overline{\partial}$ operator 
is described by the singular current that has a 
support on the wall (of the wall crossing of Lagrangian Floer cohomology). 
The wall becomes more and more dense when we 
consider the deformation of higher and higher degree.
Also in \cite{Fuk05I} the 
Feynman diagram of 0-loop is used to describe the 
scattering of the walls.
\par
After \cite{Fuk05I} had been submitted for publication, 
 Kontsevich-Soibelman \cite{KS06}
discussed the same phenomenon using 
\v Cech cohomology and studied 
the deformation of the complex structure 
as the deformation of the coordinate change 
instead of the deformation of  $\overline{\partial}$ operator.
They introduced a nilpotent group of coordinate changes 
which has a filtration by the order of the deformations.
(This group coincides with the group of $A_{\infty}$ automorphisms 
of  $H(T^n;\Lambda_0)$ which preserves the volume element.)
\cite{KS06} study the case of K3 surface mainly.
\par
After these works had been done,  
those pictures were used by Gross-Siebert, who 
had been working on a similar problem 
independently (in \cite{GS03} for example).
They proved a reconstruction theorem from 
toric degeneration in general dimension in \cite{GS07}.
\par
Thus the program 
(due to \cite{Fuk05I} etc.) to show the 
mirror symmetry which 
asserts a relation between `symplectic geoetry' and `complex geometry'
through `Morse homotopy' (or equivalently through `tropical 
geometry') is mostly
realized for the part
`Morse homotopy' $\Rightarrow$ `complex geometry', 
(as far as the part where Lagrangian submanfolds 
or coherent sheaves are not included).
The study of the other part of the 
program is making progress\footnote{For example the case of cotangent bundle is 
discussed in \cite{FSS, NZ09}.}.
\par
We discussed in this section the study of mirror symmetry 
based on the family Floer cohomology or Strominger-Yau-Zaslow conjecture.
One of the other important study of homological mirror symmetry 
of Calabi-Yau manifold is its proof in the case of quartic surface \cite{Sei03II}
by P. Seidel.
This is based on Seidel's study 
of directed $A_{\infty}$ category 
associated to the symplectic Lefschetz pencil.
We omit it and refer \cite{Sei03II}.
(Note in \cite{Sei03II} there is a warning that the proof of
a part of the results which is used in \cite{Sei03II} 
is not written up in detail. However 
now all the necessary results are established by 
\cite{Sei06}. Therefore the proof of the main results of 
\cite{Sei03II} is now completed.)
This theory of Seidel is also important in the study of 
mirror symmetry in the non-Calabi-Yau case,  which we discuss in the 
next section.

\section{Mirror symmetry in the non-Calabi-Yau case}
So far we discussed mirror symmetry for Calabi-Yau 
manifold only.
In this last section, we discuss non-Calabi-Yau case.
The mirror of a manifold
$X$ is expected to be another manifold
$\hat X$ in the Calabi-Yau case. However in a more general situation a mirror to 
a manifold $X$ becomes a pair $(\hat X, W)$ of a manifold $\hat X$ and 
a function $W$ on it. 
The function $W$ is called the {\it Landau-Ginzburg super potential}.
We explain this point here.
\par
We first consider the case when $X$ is a symplectic manifold $(X, \omega)$.
As we explained in the last section, 
the  mirror $\hat X$  of a Calabi-Yau manifold $X$
is regarded as the moduli space of skyscraper sheaves.
The homological mirror symmetry conjectures 
that  the moduli space of skyscraper sheaves  is 
identified with a moduli space of the pair $(L, b)$ 
where $L$ is a Larangian submanifold of $X$ and 
$b$ is a flat unitary connection on it ($b \in H^1(L;\sqrt{-1}\R)$).
Using it to construct $\hat X$ is an ideal of the construction of $\hat X$.
\par
In the case of Calabi-Yau manifold the space 
$\mathcal M_{\text{weak}}(L)$ (that is defined in section 2) 
coincides with 
$\mathcal M(L)$. In other words,  the potential function 
$\frak{PO}$ is a constant function $0$\footnote{We can 
prove it by the dimension counting argument.
In the case of Calabi-Yau manifold, we consider $L$ 
with vanishing Maslov index. So  
$\frak m_k(b, \cdots, b)$ never have a 
degree of the fundamental class $[L]$, if 
$b \in H^1(L)$.}.
\par
In genearal, the function  $\frak{PO}$ may not be $0$.
We regard the mirror to $(X, \omega)$ as the 
pair $(\hat X, W)$, where 
$\hat X$ is (an irreducible component of) the moduli space of the pair $(L, b)$
where $L$ is a Lagrangian submanifold and  $b \in \mathcal M_{\text{weak}}(L)$
and $W$ is the function
$$
(L, b) \mapsto W(L, b) = \frak{PO}(b).
$$
A typical case where such a construction works  
is the case when $X$ is Fano or toric.
Below we discuss on the toric case.
(\cite{Aux09} is a good reference of the material in this section.
In \cite{Kat} some applications to the algebraic geometry are also discussed.)
Toric manifold $X$ has a $T^n$ action so that its 
non-degenerate orbits are Lagrangian submanifolds.
The orbits are the fibers of the moment map $X \to P$.
The fiber $L(u)$ of the interior point $u$ of $P$ithat is a polytope in $\R^n$j
is diffeomorphic to $T^n$.
As we mentioned in section 3,  we have
 $H^1(L(u);\Lambda_0) \subset \mathcal M_{\text{weak}}(L)$.
\par
If we take an element $b$ of $H^1(L(u);\sqrt{-1}\R)$ as bounding cochain, then
$b$ corresponds to a flat unitary connection on $L(u)$. 
The totality of such $(L(u), b)$ where $u$ is a interior point of $P$, 
is the union of the dual torus of fibers $L(u)$.
We regard $\hat X$ as such totality\footnote{In this section 
we omit the argument which is related to the choice of 
the coefficient ring. Namely the difference of $\Lambda_0$
coefficient and $\C$ coefficient. I think the universal Novikov ring $\Lambda_0$
is actually the correct choice. In that case 
$\hat X$ becomes a rigid analytic space.}.
\par
We consider the function $W=\frak{PO}$ on $\hat X$.
Then Theorem \ref{thm;toric} claims that the Floer cohomology of 
$(L, b)$ is nonzero if and only if $W$ is critical at $(L, b)$.
\par
On the other hand, as we mentioned right before  Theorem 
\ref{thm;nonempty}, the Floer cohomology
$HF((L_1, b_1), (L_2, b_2))$ is defined only if $\frak{PO}(b_1) = \frak{PO}(b_2)$.
\begin{conjecture}\label{Lfswxomp}
$\mathcal{LAG}(X)$ splits into the 
direct sum of the filtered $A_{\infty}$ categories associated to each of the 
critical values of $W$.
The critical value of $W$ coincides with the 
eigenvalue of the linear map $H(X;\Lambda_0) \to H(X;\Lambda_0)$
that is given by $x \mapsto x \cup^{Q} c_1(X)$, where
$c_1(X)$ is the first Chern class and $\cup^{Q}$ is the quantum cup product.
\end{conjecture}
The author heard this conjecture for the first time in a talk 
by M. Kontsevich at  Vienna (2006).
(He was also informed that a similar conjecture was made 
by several people independently.)
\par
In the toric case,  the nontrivial part of Conjecture \ref{Lfswxomp} is stated as follows.
Let us consider $L$ and $b\in \mathcal M_{\text{weak}}(L)$ which is 
not necessary a $T^n$ orbit.
Suppose $HF((L, b), (L, b))\ne 0$.
Then Conjecture \ref{Lfswxomp} states that  $\frak{PO}(b)$ 
is equal to a critical value of $W=\frak{PO} : \hat X \to \Lambda_0$.
(Note $W$ is defined on the set of $(L(u), b)$ where $L(u)$ is a $T^n$ orbit 
and $b \in H^1(L(u);\Lambda_0)$.)
\par
We hope to prove this statement by using the fact that 
$T^n$ orbits generates the category $\mathcal{LAG}(X)$ in certain sense.
The proof of the second half is closely related to the map (\ref{openGW}).
\par
Suppose that
$c$ is a critical point of $W$.
Then $W^{-1}(c)=\hat X_c$ has a singularity.
\begin{conjecture}\label{torichm}
The derived category of the direct sum factor $\mathcal{LAG}(X)$ that corresponds to the 
critical value $c$ is isomorphic to the 
derived category of the category of the matrix factorization of a singularity  
of $\hat X_c$.
\end{conjecture}
We refer \cite{Yoshi} for the matrix factorization.
\begin{remark}
The conjectures in this section seem to be 
easier to prove than the conjectures in the last section.
We hope to solve them in a near future.
One of the reasons whey they are easier is that 
in the toric case the complex structure of $\hat X$ has 
no instanton correction.
In case when 
$X$ is not toric but is Fano,  
the complex structure of $\hat X$ may  have instanton correction.
However the examples in \cite{Aux09}
suggest that the instanton correction is simpler in 
Fano case than Calabi-Yau case, and easier to study.
\end{remark}
The  conjectures \ref{Lfswxomp} are \ref{torichm}
are a version of homological mirror symmetry.
The conjecture which is a version of classical 
mirror symmetry can be stated as follows.
The quantum cohomology of $X$ (or the Frobenius structure 
induced by the big quantum cohomology 
 by \cite{dub})  
coincides with Saito's flat structure \cite{Sai83} associated to 
$W$.
(See also \cite{ST}.)
This statement has been discussed in Givental \cite{givental1} or Hori-Vafa \cite{hori-vafa} and
has been applied successfully. 
The fact that $W$ is the potential function of  \cite{FOOO08} was
conjectured by them also. This fact is established by \cite{cho-oh, FOOO08I, fooo09}.
\par
So far we put superpotential $W$ in the complex sides.
In fact, we defined the function $W$ by using the Floer 
theory of Lagrangian submanifold on its mirror.
When we consider $W$ in the symplectic side, 
then its mirror is a complex manifold. So it is not 
natural to study Floer theory on the mirror. 
This is the reason why I discussed the case when $W$ is on the 
complex sides,  so far.
\par
In fact,  however,  more results have been obtained already 
in the case $W$ is in the symplectic side.
\par
To the pair $(X, W)$ of symplectic manifold and a function on it, 
Seidel associated a directed $A_{\infty}$ category $\mathcal{LAG}(X, W)$, 
under certain exactness condition\footnote{We say 
$M$ is exact if there exists a differential one form 
$\theta$  such that $\omega=d\theta$.
We say a Lagrangian submanifold $L$
in it is exact if there exists a function $f$ on $L$ such that $\theta = df$ holds on $L$.}.
(Seidel mentioned that his construction is 
suggested by M. Kontsevich's talk.)
The object of the Seidel's directed $A_{\infty}$ category is 
a Lagrangian submanifold that is a vanishing cycle of $W$.
The morphism space is a version of Floer cohomology. 
They are 
defined in detail in \cite{Sei06}.
\par
The mirror of such a pair $(X, W)$
is conjectured to be a compact Fano manifold $\hat X$, 
in many cases. Namely:
\begin{conjecture}
The derived category of directed 
$A_{\infty}$ category $\mathcal{LAG}(X, W)$ associated to $(X, W)$
is equivalent to the derived category of the category of coherent sheaves of $\hat X$.
\end{conjecture}
It is known that for many of the Fano manifolds, 
the derived category of the category of its coherent sheaves 
has a distinguished generator\footnote{
Such generator is given in the case of $\C P^n$ by \cite{Bei}.}.
It is conjectured that those distinguished generator 
becomes the vanishing cycle of $W$ by the above mentioned isomorphism.
This conjecture is checked by many people. (See \cite{Sei01, AKO1, Ue06, AKO, Abo} for example.)
\par
The formulation of the mirror symmetry 
for non-Calabi-Yau case is not yet completed.
\par\medskip
The Japanese version of this article was written 
in 2009. The reference below so is restricted to those which the author 
already knew at that time.
\bibliographystyle{amsplain}

\end{document}